\documentclass[twoside,a4paper,12pt]{article}

\usepackage{mathptmx}       
\usepackage{helvet}         
\usepackage{courier}        
\usepackage{type1cm}        
\usepackage{makeidx}         
\usepackage{graphicx}        
\usepackage{multicol}        
\usepackage[bottom]{footmisc}

\makeindex             
\makeindex             

\usepackage{amsmath}

\usepackage{color}

\def\pier #1{{\color{red}#1}}  
\let\pier\relax

\usepackage{cite}
\usepackage{amsmath}
\usepackage{amsthm}
\usepackage{amssymb}

\newcommand{\dega}{{\Delta_\Gamma}}
\newcommand{\bug}{{\bar u_\Gamma}}
\newcommand{\pn}{\partial_{\bf n}}
\newcommand{\pt}{\partial_t}
\newcommand{\qal}{q^{\alpha}}
\newcommand{\qgal}{q_\Gamma^{\alpha}}
\newcommand{\pal}{p^{\alpha}}

\newcommand{\ugal}{u_\Gamma^\alpha}

\newcommand{\mal}{\mu^\alpha}
\newcommand{\ral}{\rho^\alpha}
\newcommand{\rgal}{\rho_\Gamma^\alpha}
\newcommand{\bm}{\bar\mu}
\newcommand{\rg}{\rho_\Gamma}
\newcommand{\br}{\bar\rho}
\newcommand{\brg}{\bar\rho_\Gamma}
\newcommand{\ug}{{u_\Gamma}}
\newcommand{\ugan}{u_\Gamma^{\alpha_n}}
\newcommand{\mun}{\mu^{\alpha_n}}
\newcommand{\rhon}{\rho^{\alpha_n}}
\newcommand{\rgan}{\rho_\Gamma^{\alpha_n}}

\newcommand{\xiga}{{\xi_\Gamma}}
\newcommand{\rz}{{\rm I\!R}}
\newcommand{\nz}{{\rm I\!N}}

\newcommand{\dx}{{\,{\rm d}x}}
\newcommand{\dt}{{\,{\rm d}t}}
\newcommand{\ds}{{\,{\rm d}s}}
\newcommand{\dg}{{\,{\rm d}\Gamma}}

\newcommand{\oma}{{\Omega}}
\newcommand{\tint}{{\int_0^t}}
\newcommand{{\tinto}}{{\int_0^T}}
\newcommand{{\xinto}}{{\int_\Omega}}
\newcommand{{\ginto}}{{\int_\Gamma}}
\newcommand{{\texinto}}{{\int_0^T\!\!\!\int_\Omega}}
\newcommand{{\teginto}}{{\int_0^T\!\!\!\int_\Gamma}}
\newcommand{\txinto}{{\int_0^t\!\!\int_\Omega}}
\newcommand{\tginto}{{\int_0^t\!\!\int_\Gamma}}
\newcommand{\lzo}{{L^2(\Omega)}}
\newcommand{\heins}{{H^1(\Omega)}}
\newcommand{\hzwei}{{H^2(\Omega)}}
\newcommand{\Hg}{H_\Gamma}
\newcommand{\lio}{{L^\infty(\Omega)}}

\newcommand{\Vg}{{V_\Gamma}}
\newcommand{\Lg}{L^2(\Gamma)}
\newcommand{\heinsg}{H^1(\Gamma)}
\newcommand{\hzweig}{{H^2(\Gamma)}}
\newcommand{\qlzo}{{L^2(Q)}}
\newcommand{\qlio}{{L^\infty(Q)}}
\newcommand{\glzsig}{{L^2(\Sigma)}}
\newcommand{\vp}{\varphi}
\newcommand{\pig}{\pi_\Gamma}
\newcommand{\pigs}{\pi_\Gamma'}
\newcommand{\glisig}{{L^\infty(\Sigma)}}
\newcommand{\uad}{{\cal U}_{\rm ad}}
\newcommand{\overo}{\overline \Omega}

\newcommand{\CQ}{C^0(\overline Q)}
\newcommand{\CS}{C^0(\overline \Sigma)}
\newcommand{\SO}{{\cal S}_0}

\newcommand{\SALN}{{\cal S}_{\alpha_n}}
\renewcommand{\qed}{\hfill\colorbox{black}{\hspace{-0.01cm}}}

\begin{document}
\title{Optimal boundary control of a nonstandard 
Cahn--Hilliard system with dynamic boundary 
condition and double obstacle inclusions}

\author{}
\date{}
\maketitle

\begin{center}
\vspace{-15mm}
{{\sc Pierluigi Colli}\footnote{Dipartimento di Matematica  ``F. Casorati'',
Universit\`a di Pavia, \pier{Via Ferrata 5,} 27100 Pavia, Italy
(e-mail: pierluigi.colli@unipv.it) }
{\sc and J\"urgen Sprekels}\footnote{Weierstrass Institute for 
Applied Analysis and Stochastics,
Mohrenstra\ss e 39, 10117 Berlin and Department of Mathematics, 
Humboldt-Universit\"at zu Berlin, Unter den Linden 6, 10099 Berlin, \pier{Germany}
(e-mail: juergen.sprekels@wias-berlin.de) }}\\[6mm]
{{\em Dedicated to our friend Prof. Dr. Gianni Gilardi\\[1mm]
on the occasion of his 70th birthday}}\\[6mm]
{\small {\bf Key words:} optimal control; parabolic obstacle problems;\\ 
MPECs; dynamic boundary conditions; optimality conditions.\\[2mm]
}

\end{center}

\begin{abstract} 
\noindent
{\small 
In this paper, we study an optimal boundary control problem for a model for phase 
separation taking place in a
spatial domain that was introduced
\pier{by P.~Podio-Guidugli in Ric.\ Mat.\ {\bf 55} (2006), pp.~105--118}. The model consists of a strongly coupled system of nonlinear
parabolic differential inclusions, in  which products between  
the unknown functions and their time derivatives occur that are difficult to handle
analytically\pier{; the system is} complemented by initial and boundary conditions. 
For the order
parameter of the phase separation process, a dynamic boundary condition involving
the Laplace--Beltrami operator is assumed, which 
 models an additional nonconserving phase transition occurring on the surface of the domain.
We complement in this paper results that were established in the recent contribution
\pier{appeared in Evol. Equ. Control Theory {\bf 6} (2017), pp.~35--58, by the two authors and Gianni Gilardi}.  In contrast to that paper, in which differentiable potentials of 
logarithmic type were considered, we investigate here the (more difficult) case of
nondifferentiable potentials of double obstacle type. For such nonlinearities, the standard
techniques of optimal control theory to establish the existence of Lagrange multipliers
for the state constraints are known to fail.  To overcome  these difficulties, 
we employ the following line of approach: we use the results 
\pier{contained in the preprint arXiv:1609.07046 [math.AP] (2016), pp.~1--30,} 
for the case of (differentiable) logarithmic potentials and perform 
a so-called ``deep quench limit''. Using compactness and monotonicity arguments, it is shown that this strategy leads to the desired first-order necessary optimality conditions for the case of 
(nondifferentiable) double obstacle potentials.} 
\end{abstract}

\pagestyle{myheadings}
\newcommand\testopari{\sc\small Pierluigi Colli and J\"urgen Sprekels}
\newcommand\testodispari{\sc\small Boundary control of a double obstacle Cahn--Hilliard inclusion}
\markboth{\testopari}{\testodispari}
\thispagestyle{empty}
\parindent=0pt

\section{Introduction}
\noindent 
Let $\oma\subset\rz^3$ denote some open, connected and bounded
domain with smooth boundary $\Gamma$ (we should at least have $\Gamma\in C^2$), and let $T>0$ be a
fixed final time and $Q:=\oma\times (0,T)$, $\Sigma:=\Gamma\times (0,T)$.  We denote by 
$\pn$, $\nabla_\Gamma$, $\dega$, the outward normal derivative, the tangential
gradient, and the Laplace--Beltrami operator on $\Gamma$, in this order.
We study in this paper the following optimal boundary control problem:

\vspace{5mm}
(${\cal P}_0$) \quad Minimize the cost functional
\begin{align}
\label{cost}
{\cal J}((\mu,\rho,\rg),\ug):&=\frac {\beta_1} 2 \,\|\mu-\hat\mu_Q\|^2_{\qlzo}\,+\,
\frac {\beta_2} 2\,\|\rho-\hat\rho_Q\|^2_{\qlzo}\nonumber\\[1mm]
&\quad +\,\frac{\beta_3}2\,\|\rg-\hat\rho_\Sigma
\|^2_{L^2(\Sigma)}+\,\frac{\beta_4}2\,\|\rho(T)-\hat\rho_\oma\|^2_{\lzo}
\nonumber\\[1mm]
&\quad+\,\frac{\beta_5}2\,\|\rg(T)-\hat\rho_\Gamma\|^2_{L^2(\Gamma)}\,+\,
\frac{\beta_6}2\,\|\ug\|^2_{L^2(\Sigma)}
\end{align}
over a suitable set $\,\uad\subset (H^1(0,T;L^2(\Gamma))\cap \glisig)\,$ of admissible
controls $\ug$ (to be specified later), subject to the state system
\begin{align}
\label{ss1}
&(1+2g(\rho))\,\mu_t+\mu\,g'(\rho)\,\rho_t-\Delta\mu=0\quad\mbox{a.\,e. in }\,Q,\\[1mm]
\label{ss2}
&\pn\mu =0\quad\mbox{a.\,e. on }\, \Sigma,\quad \mu(0)=\mu_0 \quad\mbox{a.\,e. in }\,\oma,\\[1mm]
\label{ss3}
&\rho_t-\Delta\rho+\xi+\pi(\rho)=\mu\,g'(\rho)\quad\mbox{a.\,e. in }\,Q,\\[1mm]
\label{ss4}
&\xi\in\partial I_{[-1,1]}(\rho)\quad\mbox{a.\,e. in }\,Q,\\[1mm]
\label{ss5}
&\pn\rho +\pt\rg-\dega\rg+\xiga+\pi_\Gamma(\rg)=\ug,\quad \rg=\rho_{|\Sigma},
\quad\mbox{a.\,e. on }\,\Sigma,\\[1mm]
\label{ss6}
&\xiga\in \partial I_{[-1,1]}(\rg)\quad\mbox{a.\,e. on }\,\Sigma,\\[1mm]
\label{ss7}
&\rho(0)=\rho_0 \quad\mbox{a.\,e. in }\,\oma, \quad \rg(0)=\rho_{0_\Gamma}\quad
\mbox{a.\,e. on }\,\Gamma.
\end{align}

Here, $\beta_i$, $1\le i\le 6$, are nonnegative weights, and $\,\hat\mu_Q, \hat\rho_Q\in L^2(Q)$,
$\,\hat\rho_\Sigma\in L^2(\Sigma)$, $\,\hat\rho_\Omega\in L^2(\oma)$, and $\,\hat\rho_\Gamma\in L^2(\Gamma)$\, are prescribed target functions.

\vspace{2mm}\quad
The physical background behind the control problem (${\cal P}_0$) is the following: the state system \eqref{ss1}--\eqref{ss7} constitutes a model for phase separation taking place in the container $\oma$ and originally introduced in \cite{PG}. In this
connection, the unknowns $\,\mu\,$ and $\,\rho\,$ denote the associated \emph{chemical potential}, which 
in this particular model has to be nonnegative, and the \emph{order 
parameter} of the phase separation process, which is usually the volumetric density of one of the involved phases.
We assume that $\,\rho\,$ is normalized in such a way as to attain its values in the interval $[-1,1]$. 
The nonlinearities $\pi,\pi_\Gamma,g$ are assumed to be smooth in $[-1,1]$, and 
$\,\partial I_{[-1,1]}\,$ denotes the subdifferential of the indicator function of the
interval $[-1,1]$. As is well known, we have that
\begin{equation}
\label{indi}
I_{[-1,1]}(\rho)=\left\{\begin{array}{ll}
0&\mbox{if }\,\rho\in [-1,1]\\[1mm]
+\infty&\mbox{otherwise}
\end{array}\right.,\quad
\partial I_{[-1,1]}(\rho)=\left\{\begin{array}{ll}
(-\infty,0]&\mbox{if }\,\rho=-1\\[1mm]
\{0\}&\mbox{if }\,-1<\rho<1\\[1mm]
[0,+\infty)&\mbox{if }\,\rho=1
\end{array}
\right. .
\end{equation}
\quad The state system (\ref{ss1})--(\ref{ss7}) is singular, with highly nonlinear
and nonstandard couplings.     
It has been the subject of intensive study over the past years for the case that \eqref{ss5} is 
replaced by a zero Neumann condition. In this conncetion, we refer the reader to  
\cite{CGKPS,CGKS1,CGKS2,CGPS3, CGPS6, CGPS7, CGPS4,CGPS5}. In \cite{CGPSco}, an associated control
problem with a distributed control  in  \eqref{ss1} was investigated for the special case $g(\rho)=\rho$, 
and in~\cite{CGSco1}, the corresponding case of a boundary control for $\,\mu\,$ was studied. A nonlocal
version, in which the Laplacian $\,-\Delta\rho\,$ in \eqref{ss3} was replaced by a nonlocal operator, was
discussed in the  contributions \cite{CGS3, CGS4,CGS4neu}. 

\vspace{2mm}\quad
In all of the works cited above a zero Neumann condition was assumed for the order parameter $\,\rho$.
In contrast to this, we study in this paper the case of the dynamic boundary condition \eqref{ss5}. 
It models a nonconserving phase transition taking place on the boundary, which could be induced
by, e.\,g., an interaction between bulk and wall. The associated total free energy of the phase separation process is  
the sum of a bulk and a surface contribution and has the form
\begin{align}\label{ftot}
&{\cal F}_{\rm tot}[\pier{{}\mu(t),{}}\rho(t),\rg(t)] \nonumber\\[1mm] 
&:=\xinto \Big(I_{[-1,1]}(\rho(x,t))+\hat\pi(\rho(x,t))\,\pier{{}-\,\mu(x,t)\,g(\rho(x,t))}\,+\,\frac 12|\nabla\rho(x,t)|^2
\Big)\dx\nonumber\\[1mm] 
&+\ginto\Big(I_{[-1,1]}(\rg(x,t))+\hat\pi_\Gamma(\rg(x,t))\,\pier{{}-u_\Gamma(x,t)\,\rg(x,t)}\,+\,\frac 12|\nabla_\Gamma\rg(x,t)|^2
\Big)\dg\,,
\end{align}
for $t\in [0,T]$, where $\hat\pi(r)=\int_0^r\pi(\xi){\rm d}\xi\,$ and \,$\hat\pi_\Gamma
(r)=\int_0^r\pig(\xi){\rm d}\xi$. 
\quad
In the recent contribution \cite{CGSneu}, the state system \eqref{ss1}--\eqref{ss7} was studied
systematically concerning existence, uniqueness, and regularity. A boundary control problem resembling (${\cal P}_0$)
was solved in \cite{CGSneu2} for the case of potentials of logarithmic type.

\vspace{2mm}\quad
The mathematical literature on control problems for phase field systems involving equations
of viscous or nonviscous Cahn--Hilliard type is still scarce and quite recent. We refer in this connection to the works \cite{CFGS1, CFGS2,CGS1,CGS2,HW,wn99}. Control problems
for convective Cahn--Hilliard systems were studied in \cite{RS,ZL1,ZL2}, and a few
analytical contributions were made to the coupled Cahn--Hilliard/Navier--Stokes system
(cf. \cite{FRS,HW3, HW1,HW2}). The contribution \cite{CGRS}
dealt with the optimal control of a  Cahn--Hilliard type system arising in the modeling of
solid tumor growth. For the optimal control of Allen--Cahn equations with dynamic
boundary condition, we refer to \cite{CFS,CS} (see also \cite{Calcol}).

\vspace{2mm}\quad In this paper, we aim to employ the results established in \cite{CGSneu2}
to treat the non\-differentiable double obstacle case when $\xi,\xi_\Gamma$ satisfy the 
inclusions (\ref{ss4}), \eqref{ss6}. Our approach is guided by a strategy that was  introduced  
\pier{in \cite{CFS} by the present authors and M.H. Farshbaf-Shaker: in fact,}
we aim to derive first-order necessary
optimality conditions for the double obstacle case by 
performing a so-called ``deep quench limit'' in a family of optimal control problems
with differentiable logarithmic nonlinearities that  was treated in \cite{CGSneu2},
and for which the corresponding state systems were analyzed in \cite{CGSneu}.
The general idea is briefly explained as follows: we replace the inclusions (\ref{ss4}) and
\eqref{ss6} by the identities
\begin{equation}
\label{ss4neu}
\xi=\vp(\alpha)\,h'(\rho), \quad \xi_\Gamma=\vp(\alpha)\,h'(\rg),
\end{equation}
where $h$ is defined by
\begin{equation}
\label{defh}
h(\rho):=\left\{\begin{array}{ll}
(1-\rho)\,\ln(1-\rho)+(1+\rho)\,\ln(1+\rho)&\mbox{if } \,\rho\in (-1,1)\\[1mm]
2\,\ln(2)&\mbox{if }\,\rho\in\{-1,1\}
\end{array}\right.,
\end{equation}
and where 
 $\varphi$ is continuous and positive on $(0,1]$ and satisfies
\begin{equation}
\label{phiat0}
\lim_{\alpha\searrow 0}\,\varphi(\alpha)=0.
\end{equation}
We remark that we can simply choose $\,\varphi(\alpha)=\alpha^p\,$ for some $\,p>0$.
Now observe that $h'(y)=\ln\left(\frac{1+y}{1-y}\right)$ \,and\, $h''(y)=\frac 2 {1-y^2}>0$\, for 
$y\in (-1,1)$. Hence, in particular, we have
\begin{eqnarray}
\label{hphi}
&&\lim_{\alpha\searrow 0}\,\varphi(\alpha)\,h'(y)=0 \quad\mbox{for }\, -1<y<1,\nonumber\\[2mm]
&&\lim_{\alpha\searrow 0}\Bigl(\varphi(\alpha)\,\lim_{y\searrow -1}h'(y)\Bigr)\,=\,-\infty,
\quad \lim_{\alpha\searrow 0}\Bigl(\varphi(\alpha)\,\lim_{y\nearrow +1}h'(y)\Bigr)\,=\,+\infty\,.
\end{eqnarray}
We thus may regard the graph 
$\,\varphi(\alpha)\,h'\,$ as an approximation to the graph of the subdifferential
$\partial I_{[-1,1]}$. 

\vspace{2mm}\quad
Now, for any $\alpha>0$ the optimal control problem (later to be denoted by $({\cal P}_\alpha)$), which results if in $({\cal P}_0)$ the relations (\ref{ss4}), \eqref{ss6} are replaced by 
(\ref{ss4neu}), is of the type for which
in \cite{CGSneu2} the existence of optimal controls $u_\Gamma^\alpha\in\uad$ as well as first-order necessary optimality conditions have been derived. Proving a priori estimates (uniform in $\alpha>0$), and 
employing compactness and monotonicity arguments, we will be able to show the following existence and approximation result: whenever $\,\{u_\Gamma^{\alpha_n}\}\subset\uad$ is a sequence of optimal controls for $({\cal P}_{\alpha_n})$, where $\alpha_n\searrow 0$ as $n\to\infty$, then there exist
a subsequence of $\{\alpha_n\}$, which is again indexed by $n$, and an optimal control 
$\bar u_\Gamma\in\uad$ of
$({\cal P}_0)$ such that
\begin{equation}
\label{eq:1.15}
u_\Gamma^{\alpha_n}\to\bar u_\Gamma \quad\mbox{weakly-star in ${\cal X}$ as }\,
n\to\infty,
\end{equation}
where, here and in the following,
\begin{equation}
\label{defX} 
{\cal X}:=H^1(0,T;H_\Gamma)\cap L^\infty(\Sigma)
\end{equation}
will always denote the control space.
In other words, optimal controls for $({\cal P}_\alpha)$ are for small $\alpha>0$ likely to be `close' to 
optimal controls for $({\cal P}_0)$. It is natural to ask if the reverse holds, i.\,e., whether every optimal control for
 $({\cal P}_0)$ can be approximated by a sequence $\,\{u_\Gamma^{\alpha_n}\}\,$ of optimal controls
for $({\cal P}_{\alpha_n})$, for some sequence $\alpha_n\searrow 0$. 

\vspace{2mm}\quad
Unfortunately, we will not be able to prove such a `global' result that applies to all optimal controls for
(${\cal P}_0$). However,  a `local' result can be established. To this end, let $\bug\in\uad$ be any optimal control
for $({\cal P}_0)$. We introduce the `adapted' cost functional
\begin{equation}
\label{adcost}
\widetilde{\cal J}((\mu,\rho,\rg),\ug) \,:=\,{\cal J}((\mu,\rho,\rg),\ug)\,+\,\frac 1 2\,
\|u_\Gamma-\bar u_\Gamma\|^2_{L^2(\Sigma)}
\end{equation}
and consider for every $\alpha\in (0,1]$ the {\em adapted control problem} of minimizing $\,\widetilde{\cal J}\,$ subject to $\ug\in\uad$ and to the constraint that $(\mu,\rho,\rg)$ solves the approximating system (\ref{ss1})--(\ref{ss3}), (\ref{ss5}), \eqref{ss7}, 
(\ref{ss4neu}). It will then turn out that the following is true: 

\vspace{2mm}
(i) \,There are some sequence $\,\alpha_n\searrow 0\,$ and minimizers 
$\,{\bar u_\Gamma^{\alpha_n}}\in\uad$ of the adapted control problem 
associated with $\alpha_n$, $n\in\nz$,
such that
\begin{equation}
\label{eq:1.18}
{\bar u_\Gamma^{\alpha_n}}\to\bug\quad\mbox{strongly in $\glzsig$
as }\,n\to \infty.
\end{equation}
(ii) It is possible to pass to the limit as $\alpha\searrow 0$ in the first-order necessary
optimality conditions corresponding to the adapted control problems associated with $\alpha\in (0,1]$ in order to derive first-order necessary optimality conditions for problem $({\cal P}_0)$.

\vspace{2mm}\quad
The paper is organized as follows: in Section~2, we give a precise statement of the problem
under investigation, and we derive some results concerning the state system 
(\ref{ss1})--(\ref{ss7}) and 
its $\alpha\,$--\,approximation which is obtained if in $({\cal P}_0)$ the relations
 (\ref{ss4}) and \eqref{ss6} are replaced by the relations (\ref{ss4neu}).
In Section~3, we then prove the existence of optimal controls and the approximation 
result formulated above in
(i). The final Section~4 is devoted to the derivation of the first-order necessary 
optimality conditions, where the  strategy outlined in (ii) is employed. 

\vspace{2mm}\quad
During the course of this analysis, we will make 
repeated use of H\"older's inequality, of the elementary Young's inequality
\begin{equation}
\label{Young}
a\,b\,\le\,\gamma |a|^2\,+\,\frac 1{4\gamma}\,|b|^2\quad\forall\,a,b\in\rz \quad\forall\,\gamma>0,
\end{equation}
and of the continuity of the embeddings $H^1(\oma)\subset L^p(\oma)$ for $1\le p\le 6$.
We will also use the denotations
\begin{equation}
\label{defQt}
Q_t:=\oma\times (0,t),\quad \Sigma_t:=\Gamma\times (0,t),\quad\mbox{for }\,0<t\le T.
\end{equation}

\quad Throughout the paper, \pier{for a Banach space $\,X\,$ we denote} by $\,\|\,\cdot\,\|_X\,$ its norm and by $\,X^*\,$ its dual space. The only exemption from this rule are the norms of the
$\,L^p\,$ spaces and of their powers, which we often denote by $\|\,\cdot\,\|_p$, for
$\,1\le p\le +\infty$. By $\,\langle v,w\rangle_X\,$ we will denote
the dual pairing between elements $\,v\in X^*\,$ and $\,w\in X$. \pier{About 
the time derivative of a time-dependent function $v$, we warn the reader 
that we may use both the notation  $\pt v$ and the shorter 
one $v_t$.}

\section{General assumptions and state equations}\label{state}
\setcounter{equation}{0}
In this section, we formulate the general assumptions of the paper, and we state some preparatory results for the state system (\ref{ss1})--(\ref{ss7}) and its $\alpha$--\,approximations. To begin with,
we introduce some denotations. We set
\begin{align*}
&H:=\lzo,\quad V:=\heins,\quad W:=\{w\in\hzwei:\pn w=0\,\mbox{ on }\,\Gamma\},\\[1mm]
&\Hg:=\Lg, \quad \Vg:=\heinsg, \quad{\cal V}:=\{(v,v_\Gamma)\in V\times \Vg: v_\Gamma
=v_{|\Gamma}\},
\end{align*}
and endow these spaces with their standard norms. Notice that we have $V\subset H\subset V'$ and
$\Vg\subset \Hg\subset V_\Gamma'$, with dense, continuous and compact embeddings.

\vspace{2mm}\quad
We make the following general assumptions:
\vspace{2mm}

{\sc (A1)} \quad\,$\mu_0\in W$, \,\,$\mu_0\ge 0$ in $\overo$, \quad$\rho_0\in\hzwei$,\quad
$\rho_{0_\Gamma}:=\rho_{0_{|\Gamma}}\in \hzweig$, and
\begin{equation}\label{init}
-1\,<\,{\rm min}_{x\in\overo}\,\rho_0(x),\quad {\rm max}_{x\in\overo}\,\rho_0(x)\,<\,+1.
\end{equation}

{\sc (A2)} \quad\,$\pi,\pig\in C^2[-1,1]$; $g\in C^3[-1,1]$  is nonnegative and concave 
on $[-1,1]$.

\vspace{4mm}

{\sc (A3)} \quad\,$\uad=\left\{\ug\in {\cal X }:\,\,u_*\le \ug\le u^* \,\mbox{ a.\,e. on 
\,$\Sigma$\, and\, $\|\ug\|_{\cal X}\le R_0$}\right\},$ where 

\hspace*{12mm} $\,u_*,u^*\in \glisig\,$ and
$\,R_0>0\,$ are such that \,$\uad\not=\emptyset$.

\vspace{4mm} \quad 
Now observe that the set $\uad$ is a bounded subset of ${\cal X}$. Hence, there exists
a bounded open ball in ${\cal X}$ that contains $\uad$. For later use it is convenient to fix such a ball once and for all, noting that any other such ball could be used instead. In this sense, the following assumption is rather a denotation:

\vspace{4mm}
{\sc (A4)} \quad\,Let $\,R>0\,$ be such that $\,\uad\subset{\cal U}_R
:=\{\ug\in {\cal X}:\,\|\ug\|_{\cal X}<R\}$. 

\vspace{4mm}\quad
For the quantities entering the cost functional $\,{\cal J}\,$
(see \eqref{cost}), we assume:

\vspace{4mm}
{\sc (A5)} \quad\,The constants $\,\beta_i$, $1\le i\le 6$, are nonnegative but not
all equal to zero, 
\hspace*{12mm} and we have that $\hat\mu_Q, \hat\rho_Q\in\qlzo$, $\hat\rho_\Sigma
\in\glzsig$, $\hat\rho_\oma\in\lzo$, $\,\hat\rho_\Gamma\in\Lg$.

\vspace{5mm}
\quad We observe at this point that if {\sc (A1)}, {\sc (A2)} and $\ug\in {\cal U}_R$ hold true, then all of the 
general assumptions made in \cite{CGSneu} are satisfied provided we put, in the notation used there,
$\hat\beta=\hat\beta_\Gamma=I_{[-1,1]}$. We thus may conclude from 
\cite[Thm.~2.1 and Rem.~3.1]{CGSneu} the following well-posedness result:

\vspace{5mm}
{\sc Theorem 2.1:} \quad\emph{Suppose that the assumptions} {\sc (A1)--(A4)} {\em are
fulfilled. Then the state system} \eqref{ss1}--\eqref{ss7} {\em has for every
$\,\ug\in{\cal U}_R\,$ a unique solution $\,(\mu,\rho,\rg)\,$ with 
$\,\mu\ge 0\,$ a.\,e. in $Q$, which satisfies} 
\begin{align}
\label{regmu}
&\mu\in C^0([0,T];V)\cap L^p(0,T;W)\cap L^2(0,T;W^{2,6}(\oma))\cap \qlio
\quad\forall\,p\in [1,+\infty),\\[1mm]
\label{regmut}
&\mu_t\in L^p(0,T;H)\cap L^2(0,T;L^6(\oma))\quad\forall\,p\in [1,+\infty),\\[1mm]
\label{regrho}
&\rho\in W^{1,\infty}(0,T;H)\cap H^1(0,T;V)\cap L^\infty(0,T;H^2(\oma)),\\[1mm]
\label{regrg}
&\rg\in W^{1,\infty}(0,T;\Hg)\cap H^1(0,T;\Vg)\cap L^\infty(0,T;H^2(\Gamma)),\\[1mm]
\label{rbound}
&\rho \in [-1,1] \quad\mbox{\em a.\,e. in }\,Q, \quad \rg\in [-1,1]
\quad\mbox{\em a.\,e. on }\,\Sigma,\\[1mm]
\label{regxi}
&\xi\in L^\infty(0,T;H), \quad \xiga\in L^\infty(0,T;\Hg).
\end{align}
{\em Moreover, there is a constant $K_1^*>0$, which depends only on the data of the state
system and on $\,R$, such that}
\begin{align}\label{bounds1}
&\|\mu\|_{H^1(0,T;H)\cap C^0([0,T];V)\cap L^2(0,T;W)\cap \qlio}\,+\,
\|\rho\|_{W^{1,\infty}(0,T;H)\cap H^1(0,T;V)\cap L^\infty(0,T;H^2(\oma))}
\nonumber\\[1mm]
&+\,\|\rg\|_{W^{1,\infty}(0,T;\Hg)\cap H^1(0,T;\Vg)\cap L^\infty(0,T;H^2(\Gamma))}
\,+\,\|\xi\|_{\qlio}\,+\,\|\xiga\|_{\glisig}\,\le\,K^*_1\,,
\end{align}
{\em whenever $\,(\mu,\rho,\rg)\,$ is a solution to} \eqref{ss1}--\eqref{ss7} 
{\em which corresponds to some $\,\ug\in{\cal U}_R\,$ and satisfies} 
\eqref{regmu}--\eqref{regxi}.

\vspace{5mm}
{\sc Remark 2.2:} \quad  Thanks to Theorem 2.1, the control-to-state operator $\SO: \ug\mapsto (\mu,\rho,\rg)\,$ is
well defined as a mapping from ${\cal U}_R$ into the space specified by the regularity properties
\eqref{regmu}--\eqref{regrg}. Moreover, in view of \eqref{regrho}, it follows from well-known 
embedding results (see, e.\,g.,
\cite[Sect.~8, Cor.~4]{Simon}) that $\,\rho\in C^0([0,T];H^s(\oma))\,$ for 
$0<s<2$. In particular, we have $\,\rho\in C^0(\overline Q)$, so that $\,
\rg=\rho_{|\Gamma}\in C^0(\overline\Sigma)$.

\vspace{5mm}\quad We now turn our interest to the $\,\alpha\,$--\,approximating
system that results if we replace \eqref{ss4} and \eqref{ss6} by \eqref{ss4neu},
with $\,h\,$ given by \eqref{defh} and $\,\vp\,$ satisfying \eqref{phiat0}.  
We then obtain the following system of equations:
\begin{align}
\label{als1}
&(1+2g(\ral))\,\mu^\alpha_t+\mal\,g'(\ral)\,\rho^\alpha_t-\Delta\mal=0 \quad
\mbox{a.\,e. in }\,Q,\\[1mm]
\label{als2}
&\pn\mal=0\quad\mbox{a.\,e. on }\,\Sigma,\quad \mal(0)=\mu_0\quad\mbox{a.\,e. in }\,
\oma,\\[1mm]
\label{als3}
&\rho^\alpha_t-\Delta\ral+\vp(\alpha)\,h'(\ral)+\pi(\ral)\,=\,\mal\,g'(\ral)
\quad\mbox{a.\,e. in }\,Q,\\[1mm]
\label{als4}
&\pn\ral+\pt\rgal-\dega\rgal+\vp(\alpha)\,h'(\rgal)+\pig(\rgal)=\ugal,
\quad \rgal=\rho^\alpha_{|\Sigma} \quad\mbox{a.\,e. on }\,\Sigma,\\[1mm]
\label{als5}
&\ral(0)=\rho_0\quad\mbox{a.\,e. in }\,\oma, \quad\rgal(0)=\rho_{0_\Gamma} \quad
\mbox{a.\,e. on }\,\Gamma.
\end{align}
By virtue of \cite[Thm.~2.4]{CGSneu2}, the system \eqref{als1}--\eqref{als5} has for every $\ugal\in{\cal U}_R$
a unique solution $(\mal,\ral,\rgal)$ satisfying $\,\mal\ge 0\,$ in $Q$ and \eqref{regmu}--\eqref{regrg}.
 Moreover, there are constants $r_*(\alpha),r^*(\alpha)\in (-1,1)$, which depend
only on $R$, $\alpha$, and the data of the system, such that, for all $(x,t)\in
\overline Q$,
\begin{equation}
\label{separ}
-1<r_*(\alpha)\le\ral(x,t)\le r^*(\alpha)<1,\quad
-1<r_*(\alpha)\le\rgal(x,t)\le r^*(\alpha)<1.
\end{equation}
 Again it follows (recall Remark 2.2) that 
$\,\ral\in C^0(\overline Q)\,$ and $\,\rgal\in C^0(\overline\Sigma)$. Therefore, we may infer
from {\sc (A2)}  that there is a constant $K_2^*>0$,
which depends only on $R$ and the data of the system, such that
\begin{align}
\label{albo1}
\max_{0\le i\le 3}
\left\|g^{(i)}(\ral)\right\|_{C^0(\overline Q)}
+\,\max_{0\le i\le 2}\left(\left\|\pi^{(i)}(\ral)\right\|_{C^0(\overline Q)}\,+\,
\Bigl\|\pig^{(i)}(\rgal)\right\|_{C^0(\overline\Sigma)}\Bigr)\,\le\,K_2^*\,,
\end{align}
for every solution triple $\,(\mal,\ral,\rgal)\,$ corresponding to some $\ug\in{\cal U}_R$ and any $\alpha\in (0,1]$. Observe that a corresponding estimate cannot be concluded
for the derivatives of $\,\vp(\alpha)\,h$, since it may well happen that
$r_*(\alpha)\searrow -1$ and/or $r^*(\alpha)\nearrow +1$, as $\alpha\searrow0$.

\vspace{1mm}\quad
According to the above considerations, for every $\alpha\in (0,1]$ the solution
operator $\,{\cal S}_\alpha:\ugal\in {\cal U}_R\mapsto (\mal,\ral,\rgal)\,$ is well defined 
as a mapping into the space that is specified by the regularity properties
\eqref{regmu}--\eqref{regrg}.
We now aim to derive some
a priori estimates for $(\mal,\ral,\rgal)$ that are independent of $\alpha\in (0,1]$.
We have the following result.

\vspace{4mm}
{\sc Proposition 2.3:}\quad\,{\em Suppose that} {\sc (A1)}--{\sc (A4)} 
{\em are satisfied. Then there is some constant $K_3^*>0$, which depends only
on $R$ and on the data of the system, such
that we have: whenever $(\mal,\ral,\rgal)={\cal S}_\alpha(\ugal)$ for some 
$\ugal\in {\cal U}_R$ and some
$\alpha\in (0,1]$, then it holds that}
\begin{align}
\label{albo2}
&\|\mu^\alpha\|_{H^1(0,T;H)\cap C^0([0,T];V)\cap L^2(0,T;W)\cap \qlio}
\nonumber\\[1mm]
&+\,\|\pier{\rho^\alpha}\|_{W^{1,\infty}(0,T;H)\cap H^1([0,T];V)
\cap L^\infty(0,T;H^2(\Gamma))}\nonumber\\[1mm]
&+\,\|\rho_\Gamma^\alpha\|_{W^{1,\infty}(0,T;H_\Gamma)\cap H^1([0,T];V_\Gamma)
\cap L^\infty(0,T;H^2(\Gamma))}\nonumber\\[1mm]
&+\|\vp(\alpha)\,h'(\ral)\|_{L^\infty(0,T;H)}\,+\,\|\vp(\alpha)\,h'(\rgal)\|_{L^\infty(0,T;\Hg)}\,\le\,K_3^*\,.
\end{align} 

\vspace{2mm}
{\sc Proof:}\,\quad Let  $\ugal\in {\cal U}_R$ and $\alpha\in (0,1]$ be arbitrary and $(\mal,\ral,\rgal)={\cal S}_\alpha(\ugal)$. The result will be established in a series of a priori
estimates. To this end, we will in the following \pier{denote by $C>0$ constants} 
that may depend on the quantities mentioned in the statement, but not on 
$\alpha\in (0,1]$. For the sake of a better readability, we will omit the superscript $\,\alpha\,$ of
$(\mal,\ral,\rgal)$
 during the estimations, writing it only at the end of each
estimate. We will also make repeated use of the general bounds \eqref{albo1} without
further reference. 

\vspace{4mm}
\underline{\sc First estimate:} 

\vspace{1mm}
First, note that $\,\,\pt((\mbox{$\frac 12$}+g(\rho))\,\mu^2)\,=\,
(1+2g(\rho))\,\mu_t\,\mu\,+\,g'(\rho)\,\rho_t\,\mu^2$.
Thus, multiplying \eqref{als1} by $\mu$ and integrating over $Q_t$, where $t\in (0,T]$,
we find the estimate
\begin{align}
\label{asti1}
\xinto\left(\mbox{$\frac 12$}+g(\rho(t))\right)|\mu(t)|^2\dx
\,+\,\txinto|\nabla\mu|^2\dx\ds\,=\, 
\xinto\left(\mbox{$\frac 12$}+g(\rho_0)\right)|\mu_0|^2\dx\,.
\end{align}
Hence, as $g(\rho)\ge 0$ by {\sc (A2)}, it follows that
\begin{equation}\label{asti2}
\left\|\mal\right\|_{L^\infty(0,T;H)\cap L^2(0,T;V)}\,\le\,C\,
\quad\forall\,\alpha\in (0,1].
\end{equation}

\vspace{4mm}
\underline{\sc Second estimate:}

\vspace{1mm} Next, we multiply \eqref{als3} by $\,\vp(\alpha)\,h'(\ral)\,$ and 
integrate over $Q_t$ and by parts, where $t\in (0,T]$. We obtain the identity
\begin{align}
\label{asti3}
&\vp(\alpha)\xinto h(\rho(t))\dx\,+\,\vp(\alpha)\ginto h(\rg(t))\dg\,+\,
\txinto|\vp(\alpha)\,h'(\rho)|^2\dx\ds\nonumber\\
&{}+\tginto|\vp(\alpha)\,h'(\rg)|^2\dg\ds
\,+\,\vp(\alpha)\txinto h''(\rho)\,|\nabla\rho|^2\dx\ds\nonumber\\
&{}+\,\vp(\alpha)\tginto h''(\rg)\,|\nabla_\Gamma\rg|^2\dg\ds\nonumber\\
&{}=\,\vp(\alpha)\xinto h(\rho_0)\dx\,+\,\vp(\alpha)\ginto h(\rho_{0_\Gamma})\dg
\nonumber\\
&\quad
+\txinto(\mu\,g'(\rho)-\pi(\rho))\,\pier{\vp(\alpha)}\,h'(\rho)\dx\ds\nonumber\\
&\quad +\tginto(\ugal-\pig(\rg))\,\pier{\vp(\alpha)}\,h'(\rg)\dg\ds\,.
\end{align}
Obviously, all of the terms on the left-hand side are nonnegative, while
the first two summands on the right-hand \pier{side} are bounded independently 
of $\alpha\in (0,1]$.
Thus, applying \pier{H\"older's and Young's inequalities to the last two integrals in \eqref{asti3}}, and invoking \eqref{albo1} and \eqref{asti2}, we readily find that
\begin{equation}
\label{asti4}
\left\|\vp(\alpha)\,h'(\ral)\right\|_{\qlzo}\,+\,\left\|\vp(\alpha)\,h'(\rgal)
\right\|_{\glzsig}\,\le\,C\quad\forall\,\alpha\in (0,1]. 
\end{equation}

\vspace{3mm}
\underline{\sc Third estimate:} 

\vspace{1mm}
We now add $\,\rho\,$ on both sides of \eqref{als3} and $\,\rg\,$ on both sides
of (\ref{als4}). Then we multiply the first resulting equation by $\rho_t$ and
integrate over $Q_t$, where $t\in (0,T]$. Employing \eqref{albo1}, we then
obtain an inequality of the form
\begin{align}
\label{asti5}
&\txinto|\rho_t|^2\dx\ds\,+\,\tginto|\pt\rg|^2\dx\ds\,+\,\frac 1 2\,
\left(\|\rho(t)\|_V^2\,+\,\|\rg(t)\|_{\Vg}^2\right)\nonumber\\
&\le\,\frac 12\left(\|\rho_0\|^2_V\,+\|\rho_{0_\Gamma}\|^2_{\Vg}\right)
\,+\txinto|\rho_t|\,(|\rho|\,+\,|\vp(\alpha)\,h'(\rho)|\,+\,C(1+|\mu|))\dx\ds
\nonumber\\
&\quad \,+\tginto|\pt\rg|\,(|\rg|\,+\,|\vp(\alpha)\,h'(\rg)|\,+\,|\ugal|)\dg\ds\,.
\end{align}
Using (A1), \eqref{asti2}, and \eqref{asti4}, and employing Young's inequality
and Gronwall's lemma, we thus conclude that 
\begin{equation}\label{asti6}
\left\|\ral\right\|_{H^1(0,T;H)\cap L^\infty(0,T;V)}\,+
\,\left\|\rgal\right\|_{H^1(0,T;\Hg)\cap L^\infty(0,T;\Vg)} \,\le\,C\quad\forall 
\,\alpha\in (0,1]\,.
\end{equation}

\vspace{5mm}
\underline{\sc Fourth estimate:} 

\vspace{1mm} 
We now take advantage of the estimates \eqref{albo1}, \eqref{asti2}, \eqref{asti4} and \eqref{asti6}. Indeed,
comparison in \eqref{als3} yields that
\begin{align}
\label{asti7}
\left\|\Delta\rho\right\|_{\qlzo}\,\le\,C\,.
\end{align}   
Now observe that, owing to \cite[Thm.~3.2, p.~1.79]{brez}, we have the estimate
\begin{equation*}
\int_0^T\|\rho(t)\|^2_{H^{3/2}(\oma)}\dt\,\le\,C\int_0^T\left(\|\Delta\rho(t)\|_H^2\,+
\|\rg(t)\|_{\Vg}^2\right)\dt,
\end{equation*}
so that
\begin{equation}\label{asti8}
\left\|\rho\right\|_{L^2(0,T;H^{3/2}(\oma))}\,\le\,C.
\end{equation}
Hence, by the trace theorem (cf. \cite[Thm.~2.27, p.~1.64]{brez}), \pier{we infer that}
\begin{align}
\label{asti9}
\|\pn\rho\|_{L^2(0,T;\Hg)}\,\le\,C,
\end{align} 
whence, by comparison in \eqref{als4},
\begin{equation}
\label{asti10}
\|\dega\rg\|_{L^2(0,T;\Lg)}\,\le\,C.
\end{equation}
Thus, by the boundary version of elliptic estimates, \pier{we deduce that}
\begin{equation}\label{asti11}
\|\rg\|_{L^2(0,T;\hzweig)}\,\le \,C,
\end{equation}
whence, by virtue of standard elliptic theory\pier{, it turns out that}
\begin{equation}\label{asti12}
\|\rho\|_{L^2(0,T;\hzwei)}\,\le\,C.
\end{equation} 
Since the embeddings 
$$(H^1(0,T;H)\cap L^2(0,T;\hzwei))\subset C^0([0,T];V)$$ and $$(H^1(0,T;\Hg)
\cap L^2(0,T;\hzweig))\subset C^0([0,T];\Vg)$$
are continuous, we have thus shown the estimate
\begin{align}
\label{asti13}
\left\|\ral\right\|_{C^0([0,T];V)\cap L^2(0,T;\hzwei)}\,+\,\left\|\rgal\right\|
_{C^0([0,T];\Vg)\cap L^2(0,T;\hzweig)} \,\le\,C \quad\forall\,\alpha\in (0,1].
\end{align}

\vspace{5mm}
\underline{\sc Fifth estimate:}

\vspace{1mm}\noindent
In this step of the proof, we \pier{adopt} a formal argument that can be made rigorous by 
using finite differences in time. Namely, we differentiate \eqref{als3}
formally with respect to time, multiply the resulting identity by $\,\rho_t$, and integrate over $Q_t$,
where $0<t\le T$, and (formally) by parts. We then arrive at an inequality of the form
\begin{align}
\label{asti14}
&\frac 12\left(\|\rho_t(t)\|_H^2+\|\pt\rg(t)\|_{\Hg}^2\right)\,+\txinto|\nabla\pt\rho|^2\dx\ds
\,+\tginto|\nabla_\Gamma\pt\rg|^2\dg\ds\nonumber\\
&+\,\vp(\alpha)\txinto h''(\rho)|\rho_t|^2\dx\ds\,+\,\vp(\alpha)\tginto h''(\rg)|\pt\rg|^2\dg\ds
\nonumber\\
&\le\,\frac 12\left(\|\rho_t(0)\|_H^2+\|\pt\rg(0)\|_{\Hg}^2\right)\,+\,\sum_{j=1}^4 I_j,
\end{align}    
where the expressions $\,I_j$, $1\le j\le 4$, will be specified and estimated below. Notice
that all of the terms on the left-hand side are nonnegative. At first, 
using (A1), (A2), the trace theorem, and the fact that $\ugal\in\uad$, we find that
\begin{align}
\label{asti15}
&\|\rho_t(0)\|_H\,\le\,\|\Delta\rho_0-\vp(\alpha)\,h'(\rho_0)-\pi(\rho_0)+
\mu_0\,g'(\rho_0)\|_H\,\le\,C,\nonumber\\[1mm]
&\|\pt\rg(0)\|_H\,\le\,\|\pn\rho_0\|_H\,+\,\|\Delta_\Gamma\rho_{0_\Gamma}
-\vp(\alpha)\,h'(\rho_{0_\Gamma})-\pig(\rho_{0_\Gamma})+\ugal(0)\|_H\,\le\,C.
\end{align}
Next, recalling \eqref{albo1} and \eqref{asti6}, we have that
\begin{equation}\label{asti16}
I_1:=-\txinto \pi'(\rho)\,|\rho_t|^2\dx\ds\,\le\,C,
\end{equation}  
as well as, by also using Young's inequality,
\begin{equation}
\label{asti17}
I_4:=\tginto(\pt\ugal-\pig'(\rg)\,\pt\rg)\,\pt\rg\dg\ds\,\le\,C.
\end{equation}
In addition, since $\,\mu\,g''(\rho)\le 0$, \pier{it turns out that}
\begin{align}
\label{asti18}
&I_2:=\txinto \mu\,g''(\rho)\,|\rho_t|^2\dx\ds\,\le\,0.
\end{align} 

The estimation of the remaining term
$$
I_3:=\txinto\mu_t\,g'(\rho)\,\rho_t\dx\ds
$$
is more delicate. To this end, we use the identity \pier{(cf.~\eqref{als1})}
$$\mu_t\,=\,(1+2g(\rho))^{-1}(\Delta\mu-\mu\,g'(\rho)\,\rho_t),$$
where, obviously, $\,1/(1+2g(\rho))\le 1$. Substitution of this 
identity and integration by parts yield that
\begin{align}
I_3&= \txinto \frac 1{1+2g(\rho)} \,
  \bigl[ \Delta \mu - \mu \,g'(\rho)\, \rho_t]\,
  g'(\rho)\,\rho_t \dx\ds
  \nonumber
  \\[1mm]
&=\,-\txinto\nabla\mu(s) \cdot \nabla \Bigl(
  \frac {g'(\rho)\,\rho_t} {1+2g(\rho)} \Bigr)\dx\ds
    \,  - \txinto \frac {(g'(\rho))^2}{1+2g(\rho)} \, \mu\,|\rho_t|^2\dx \ds\,, 
  \label{asti19}
\end{align}
where the second summand on the right is obviously nonpositive.
We thus obtain the inequality
\begin{align}
\label{asti20}
 I_3\,\le \,C\txinto|\nabla\mu|\,|\nabla\rho_t|\dx\ds\,+\,C\txinto|\nabla\mu|\,
|\nabla\rho|\,|\rho_t|\dx\ds :=J_1+J_2\,.
\end{align}
Obviously, owing to Young's inequality and \eqref{asti2}, \pier{we infer that}
\begin{equation}
\label{asti21}
J_1\,\le\,
\frac 14 \,\txinto  |\nabla\rho_t|^2 \dx\ds\,+\,C.
\end{equation}
On the other hand, thanks to H\"older's and Young's inequalities, we also have that
\begin{align}
\label{asti22}
J_2  &\le \,C \int_0^t \|\nabla\mu(s)\|_2\, \|\nabla\rho(s)\|_4 \, \|\rho_t(s)\|_4 \dx\ds
  \nonumber  \\
  & \leq \,\frac 14 \tint\|\rho_t(s)\|_V^2\ds\,+\,
  C\tint\|\nabla\mu(s)\|_H^2\,\|\nabla\rho(s)\|_V^2\ds\nonumber\\
  & \leq \,C\,+\,\frac 14\txinto |\nabla\rho_t|^2\dx\ds\,
  + \,C\tint  \|\nabla\mu(s)\|_H^2 \|\nabla\rho(s)\|_V^2 \ds\,.
 \end{align}

\quad The last integral cannot be controlled in this form. We thus try to estimate
the expression $\,\|\nabla\rho(s)\|_V^2	\,$ in terms of the expressions
\,$\,\|\pt\rho(s)\|_H^2\,$ and $\,\|\pt\rg(s)\|^2_{\Hg}\,$ which can be handled 
using the first summand on the left-hand side of  \eqref{asti14}. To this end, we use the 
regularity theory for linear elliptic equations and \eqref{asti13} to deduce that 
\begin{equation}
\label{asti23}
\|\nabla\rho(s)\|_V^2\,\le\, C\left(\|\rho(s)\|_V^2\,+\,\|\Delta\rho(s)\|_H^2\right)
\,\le\,C\left(1\,+\,\|\Delta\rho(s)\|_H^2\right).
\end{equation}
We now multiply, just as in the second estimate above, \eqref{als3} by 
$\,\vp(\alpha)\,h'(\rho(s))$, but this time we only integrate over $\oma$.
We then obtain, for almost every $s\in (0,t)$,          
\begin{align}
\label{asti24}
&\|\varphi(\alpha)\,h'(\rho(s))\|_H^2\,+
\,\|\varphi(\alpha)\,h'(\rg(s))\|_{\Hg}^2\,+\,\vp(\alpha)
\xinto h''(\rho(s))\,|\nabla\rho(s)|^2\dx\nonumber\\
&+\,\vp(\alpha) \ginto h''(\rg(s))\,|\nabla_\Gamma\rg(s)|^2\dg\nonumber\\
&=\,\xinto\vp(\alpha)\,h'(\rho(s))\left(-\rho_t(s)-\pi(\rho(s))+\mu(s)\,g'(\rho(s))\right)\dx
\nonumber\\
&\quad +\ginto\vp(\alpha)\,h'(\rg(s))\left(-\pt\rg(s)-\pig(\rg(s))+\pt\ugal(s)\right)\dg\,,
\end{align}
whence, thanks to the already proven estimates and to Young's inequality,
\begin{align}
\label{asti25}     
&\|\varphi(\alpha)\,h'(\rho(s))\|^2_H\,+\,\|\varphi(\alpha)\,h'(\rg(s))\|_{\Hg}^2\,
\le\,C\left(1+\|\pt\rho(s)\|_H^2+\|\pt\rg(s)\|^2_{\Hg}\right)\nonumber\\[1mm]
&\mbox{for a.\,e. $s\in (0,t)$.}
\end{align}
Comparison in \eqref{als3} then yields that
\begin{equation}
\label{asti26}
\|\Delta\rho(s)\|_H^2\,\le\,C\left(1+\|\pt\rho(s)\|_H^2+\|\pt\rg(s)\|^2_{\Hg}\right)
\quad\mbox{for a.\,e. $s\in (0,t)$.}
\end{equation}
Combining the estimates \eqref{asti20}--\eqref{asti26}, we have thus shown that
\begin{align}
I_3\,\le\,C+\frac 12 \txinto|\nabla\rho_t|\dx\ds\,+\,C\tint\|\nabla\mu(s)\|_H^2
\left(\|\rho_t(s)\|_H^2+\|\pt\rg(s)\|_H^2\right)\dx\ds\,,
\label{asti27}
\end{align}
where the mapping $\,\,s\mapsto\|\nabla\mu(s)\|_H^2\,\,$ is known to be bounded in
$L^1(0,T)$, uniformly with respect to $\alpha\in (0,1]$. We thus may combine 
\eqref{asti14}--\eqref{asti18} with \eqref{asti27} to infer from Gronwall's lemma that
\begin{align}
\label{asti28}
\left\|\ral\right\|_{W^{1,\infty}(0,T;H)\cap H^1(0,T;V)}\,+\,
\left\|\rgal\right\|_{W^{1,\infty}(0,T;\Hg)\cap H^1(0,T;\Vg)}\,\le\,C
\quad\forall\,\alpha\in (0,1].
\end{align}
Therefore, we can conclude from \eqref{asti25} and \eqref{asti26} that also, for all $\alpha\in (0,1]$,
\begin{equation}
\label{asti29}
\|\varphi(\alpha)\,h'(\ral)\|_{L^\infty(0,T;H)}\,+\,\|\varphi(\alpha)\,h'(\rgal)\|_{L^\infty(0,T;\Hg)}\,
+\,\|\Delta\ral\|_{L^\infty(0,T;H)}\,\le\,C.
\end{equation}
Since we already know from \eqref{asti13} the bound for $\,\|\rgal\|_{C^0([0,T];\Vg)}$, 
we can follow the same chain of estimates as in the fourth a priori estimate above,
eventually obtaining that
\begin{align}
\label{asti30}
\left\|\ral\right\|_{L^\infty(0,T;\hzwei)} \,+\,\left\|\rgal\right\|_{L^\infty(0,T;\hzweig)}
\,\le\,C \quad\forall\,\alpha\in (0,1].
\end{align}

\vspace{5mm}
\underline{\sc Sixth estimate:}

\vspace{1mm}
Next, we multiply \pier{\eqref{als1}} by $\,\mu_t\,$ and integrate over $Q_t$, where $t\in (0,T]$. Recalling that
$\,g(\rho)\,$ is nonnegative, and using H\"older's and Young's inequalities, we obtain from 
(A1) that
\begin{align}
\label{asti31}
&\txinto|\mu_t|^2\dx\ds\,+\,\frac 12\,\|\nabla\mu(t)\|_H^2\,\le\,\frac 12\,\|\nabla\mu_0\|_H^2
\,+\,C\txinto|\mu_t|\,|\mu|\,|\rho_t|\dx\ds\nonumber\\
&\le\,C\,+\,C\tint\|\mu_t(s)\|_2\,\|\mu(s)\|_4\,\|\rho_t(s)\|_4\ds\nonumber\\
&\le\,C\,+\,\frac 12 \txinto|\mu_t|^2\dx\ds\,+\,C\tint\|\rho_t(s)\|_V^2\,\|\mu(s)\|_V^2\ds\,,
\end{align}
where, owing to \eqref{asti28}, the mapping $\,\,s\mapsto\|\rho_t(s)\|_V^2\,\,$ 
is bounded in $L^1(0,T)$, uniformly in $\alpha\in (0,1]$. We thus can infer 
from Gronwall's lemma that 
\begin{align}
\label{asti32}
\|\mu\|_{H^1(0,T;H)\cap L^\infty(0,T;V)}\,\le\,C.
\end{align}
Comparison in \eqref{als1} then shows that also
\begin{equation}\label{asti33}
\|\Delta\mu\|_{L^2(0,T;H)}\,\le\,C,
\end{equation}
whence, by virtue of standard elliptic estimates,
\begin{equation}\label{asti34}
\|\mu\|_{L^2(0,T;W)}\,\le\,C.
\end{equation}
Since the embedding $(H^1(0,T;H)\cap L^2(0,T;\hzwei))\subset C^0([0,T];V)$ is continuous,
we have thus shown the estimate
\begin{equation}
\label{asti35}
\left\|\mal\right\|_{H^1(0,T;H)\cap C^0([0,T];V)\cap L^2(0,T;W)}\,\le\,C
\quad\forall\,\alpha\in (0,1].
\end{equation}
\quad Next, we use \pier{the continuity of the} embedding $$(L^\infty(0,T;H)\cap L^2(0,T;V))
\subset L^{7/3}(0,T;L^{14/3}(\oma)),$$ which, in view of 
\eqref{asti28}, implies that
\begin{equation}
\label{asti36}
\left\|\rho^\alpha_t\right\|_{L^{7/3}(0,T;L^{14/3}(\oma))}\,\le\,C\quad\forall\alpha\in (0,1].
\end{equation}
With this estimate shown, we may argue as in the proof of \cite[Thm.~2.3]{CGPS3} to
conclude that 
\begin{equation}
\label{asti37}
\left\|\mal\right\|_{L^\infty(Q)}\,\le\,C\quad\forall\,\alpha\in (0,1].
\end{equation}

\pier{Hence}, the assertion is completely proved.\qed

\section{Existence and approximation of optimal controls}
\setcounter{equation}{0}

In this section, we aim to approximate optimal pairs of (${\cal P}_0$). To this end, we consider for $\alpha\in (0,1]$ the optimal control problem

\vspace{5mm}
(${\cal P}_\alpha$) \quad Minimize the cost functional
$\,{\cal J}((\mal,\ral,\rgal),\ugal)\,$ for $\,\ugal\in\uad$,  subject 
\hspace*{13mm}to the state system \eqref{als1}--\eqref{als5}.

\vspace{5mm}
According to \cite[Thm.~4.1]{CGSneu2}, this optimal control problem has an
optimal pair $((\mal,\ral,\rgal),$ $\ugal)$, for every $\alpha\in (0,1]$.
Our first aim in this section is to prove the following approximation result:

\vspace{5mm}
{\sc Theorem 3.1:}\quad\,{\em Suppose that the assumptions} {\sc (A1)}--{\sc (A5)} 
{\em are satisfied, and let \pier{the} sequences $\{\alpha_n\}\subset (0,1]$ and 
$\{\ugan\}\subset \uad$ be given such that $\alpha_n\searrow 0$ and
$\ugan\to\ug$ weakly-star in ${\cal X}$ for some $\ug\in\uad$. Then it holds,
for $\,(\mun,\rhon,\rgan)={\cal S}_{\alpha_n}(\ugan)$, $n\in\nz$,} 
\begin{align}
\label{conmu}
\mun&\to \mu\quad\mbox{\em weakly-star in }\,
\pier{H^1(0,T;H)\cap L^\infty(0,T;V)\cap L^2(0,T;W)\cap L^\infty(Q)} , \\[1mm]
\label{conrho}
\rhon&\to \rho\quad\mbox{\em weakly-star in }\,W^{1,\infty}(0,T;H)\cap H^1(0,T;V)
\cap L^\infty(0,T;\hzwei),\\[1mm]
\label{conrg}
\rgan&\to \rg\quad\mbox{\em weakly-star in }\,W^{1,\infty}(0,T;{\Hg})\cap H^1(0,T;{\Vg})
\cap L^\infty(0,T;\hzweig),
\end{align}
\pier{\em as well as} 
\begin{align}
\label{conxi}
\vp(\alpha_n)\,h'(\rhon)&\to\xi\quad\mbox{\em weakly-star in }\,L^\infty(0,T;H),\\[1mm]
\label{conzig}
\vp(\alpha_n)\,h'(\rgan)&\to\xi_\Gamma\quad\mbox{\em weakly-star in }\,L^\infty(0,T;\Hg),
\end{align}
{\em where $(\mu, \rho,\rg,\xi,\xi_\Gamma)$ is the unique solution to the state
system} \eqref{ss1}--\eqref{ss7} {\em associated with $\ug$. Moreover, with
$\,{\cal S}_0(\ug)=(\mu,\rho,\rg)\,$ it holds that}
\begin{align}
\label{conj1}
{\cal J}({\cal S}_0(\ug),\ug)\,&\le\,\liminf_{n\to\infty}\,{\cal J}({\cal S}_{\alpha_n}
(\ugan),\ugan),\\[1mm]
\label{conj2}
{\cal J}({\cal S}_0(v_\Gamma),v_\Gamma)\,&=\lim_{n\to\infty}\,{\cal J}({\cal S}_
{\alpha_n}(v_\Gamma),v_\Gamma)\quad\forall\,v_\Gamma\in\uad.
\end{align}
 
\vspace{3mm}
{\sc Proof:} \,\quad
Let $\,\{\alpha_n\}\subset (0,1]\,$ be any sequence such that $\alpha_n\searrow 0$ as $n\to\infty$, and suppose that $\ugan\to\ug$ weakly-star in ${\cal X}$ for some $\ug\in\uad$. By virtue of Proposition 2.3, there are a subsequence of $\{\alpha_n\}$, which is 
again indexed by $n$, and some quintuple $\,(\mu,\rho,\rg,\xi,\xi_\Gamma)\,$ such that
the convergence results \eqref{conmu}--\eqref{conzig} hold true. In particular, we 
have $\,\mu(0)=\mu_0\,$ and $\,\rho(0)=\rho_0$. Moreover, from standard compact 
embedding results (cf. \cite[Sect.~8, Cor.~4]{Simon}) we can infer that
\begin{align}
\label{strcon1}
\mun&\to\mu\quad\mbox{strongly in }\,C^0(0,T;H)\cap L^2(0,T;V),\\[1mm]
\label{strcon2}
\rhon&\to\rho\quad\mbox{strongly in }\,\CQ,
\end{align}
\pier{also including}
\begin{equation}
\label{strcon3}
\rgan\to\rg\quad\mbox{strongly in }\,\CS,
\end{equation}
whence we infer that $\,\rg=\rho_{|\Sigma}$. \pier{Therefore}, we obviously have that
\begin{align}
\label{strcon4}
\Psi(\rhon)&\to\Psi(\rho)\quad\mbox{strongly in $\,\CQ$, for }\,\Psi\in\{g,g',\pi\},\\[1mm]
\label{strcon5}
\pig(\rgan)&\to\pig(\rg)\quad\mbox{strongly in $\,\CS$,}
\end{align} 
and \pier{\eqref{conrho}} implies that $\,\pn\rhon\to\pn\rho\,$ weakly in $L^2(\Sigma)$.
\pier{Further}, we easily verify that, at least weakly in $L^1(Q)$,
\begin{align}
&g(\rhon)\,\mu_t^{\alpha_n}\to g(\rho)\,\mu_t,\quad \mun\,g'(\rhon)\,\rho_t^{\alpha_n}\to
\mu\,g'(\rho)\,\rho_t,\quad \mun\,g'(\rhon)\to\mu\,g'(\rho).
\end{align}

\vspace{2mm}
\quad Combining the above convergence results, we may pass to the limit
as $n\to\infty$ in the equations (\ref{als1})--(\ref{als5}) (written for 
$\alpha=\alpha_n$) to find that the
quintuple $(\mu,\rho, \rg, \xi,\xiga)$ satisfies the equations \eqref{ss1}--(\ref{ss3}), 
\eqref{ss5}, and \eqref{ss7}. In addition, we have $\mu\ge 0$ in $Q$, and the 
properties \pier{in (\ref{rbound}) are fulfilled. We also notice that the regularities 
in \eqref{regmu}--\eqref{regmut} follow from $\mu_0\in W$ (cf.~(A1)) and the regularity 
theory for solutions to linear uniformly parabolic equations with continuous coefficients
and right-hand side in $L^\infty(0,T;H)\cap L^2(0,T;L^6(\oma))$ (comments are given 
in~\cite[Section~3, Step~4 and Remark~3.1]{CGSneu}). Then, in} order to show that 
the quintuple $(\mu,\rho,\rg,\xi,\xi_\Gamma)$ is in fact the unique solution to 
problem (\ref{ss1})--(\ref{ss7}) corresponding to $\ug$, 
it remains to
 show that $\,\xi\in\partial I_{[-1,1]}(\rho)\,$ a.\,e. in $Q$
and $\,\xiga\in \partial I_{[-1,1]}(\rg)\, $ a.\,e. in $\Sigma$. 

\vspace{2mm}
Now, recall that $\,h\,$ is convex \pier{in $[-1,1]$ and both $\,h\,$ and $\,\vp\,$ are} nonnegative.
We thus have, for every $n\in\nz$,
\begin{align}
\label{t3.1}
\pier{0}\ &\pier{\leq 
\pier{\texinto}\varphi(\alpha_n)\,h(\rho^{\alpha_n})\dx\,\dt}\nonumber\\[1mm]
&\le\pier{\texinto}\varphi(\alpha_n)\,\pier{h(z)}\dx\,\dt\,+
\,\pier{\texinto}\varphi(\alpha_n)\,h'( \rho^{\alpha_n})\, (\pier{\rho^{\alpha_n} - z})
\dx\,\dt\nonumber\\[1mm]
&
\hskip2cm \pier{\mbox{for all }\,z\in {\cal K}:=\{v\in {L^2(Q)}:|v|\leq 1\text{ a.e. in }Q\}\,.} 
\end{align}

Thanks to (\ref{phiat0}), the \pier{first integral on the central line of \eqref{t3.1} tends to zero as $n\to\infty$.} Hence, invoking  
\pier{(\ref{conxi}) and (\ref{strcon2})},  the passage to the limit as $n\to\infty$ yields
\begin{equation}
\label{t3.2}
\pier{\texinto}\xi\,( \pier{\rho} -z)\,\dx\,\dt\,\geq 0\quad\forall z\in {\mathcal{K}} . 
\end{equation}
Inequality (\ref{t3.2}) entails that $\xi$ is an element of the subdifferential of the extension $\mathcal{I} $ of $ I_{[-1,1]}$ to $L^2(Q)$, which means that $\xi \in \partial \mathcal{I} (\rho)$ or, equivalently (cf.~\cite[Ex.~2.3.3., p.~25]{Brezis}),  
$\xi\in\partial I_{[-1,1]}(\rho)$ {a.\,e. in $Q$}. Similarly, we can prove that 
$\xi_\Gamma\in\partial I_{[-1,1]}(\rho_\Gamma)$ a.\,e. in $\Sigma$.

\vspace{1mm}
\quad We have thus shown that, for a suitable subsequence of $\{\alpha_n\}$,
we have the convergence properties \eqref{conmu}--\eqref{conzig}, where $(\mu,\rho,
\rg,\xi, \xi_\Gamma)$ is a solution to the state system \eqref{ss1}--\eqref{ss7}.
But this solution is known to be unique, which entails that the above convergence
properties are valid for the entire sequence. This finishes the proof of the
first claim of the theorem.

\vspace{2mm}
\quad It remains to show the validity of \eqref{conj1} and \eqref{conj2}. 
In view of \eqref{conmu}--\eqref{conrg}, the inequality \eqref{conj1}  is an immediate consequence of the 
weak sequential semicontinuity properties of the cost functional ${\cal J}$. To establish
the identity \eqref{conj2}, let $v_\Gamma\in \uad$ be arbitrary  and  put
$(\mu^{\alpha^n},\rho^{\alpha_n},\rg^{\alpha_n})=\SALN(v_\Gamma)$, for $n\in\nz$. 
Taking Proposition 2.3 into account, and arguing as in the first part of this proof,
we can conclude that $\SALN(v_\Gamma)$ converges to $(\mu,\rho,\rg)=\SO(v_\Gamma)$ in the sense
of \eqref{conmu}--\eqref{conrg} \pier{and \eqref{strcon1}--\eqref{strcon3}}. In particular, we have
$$\SALN(v_\Gamma)\to \SO(v_\Gamma)\quad\mbox{strongly in }\,C^0([0,T];H)\times C^0([0,T];H)
\times C^0([0,T];\Hg).
$$
As the cost functional ${\cal J}$ is obviously continuous in the variables $(\mu,\rho,\rg)$
with respect to the strong topology of $\,C^0([0,T];H)\times C^0([0,T];H)\times C^0([0,T];\Hg)$,
we may thus infer that \eqref{conj2} is valid. \qed

\vspace{7mm}
{\sc Corollary 3.2:} \quad\,{\em The optimal control problem} (${\cal P}_0$) {\em has a least one solution.}

\vspace{3mm}
{\sc Proof:} \quad Pick an arbitrary sequence $\{\alpha_n\}$ such that $\alpha_n\searrow0$ as $n\to\infty$.  
Then, by virtue of \cite[Thm.~4.1]{CGSneu2}, the optimal control problem (${\cal P}_{\alpha_n}$) has for
every $n\in\nz$ an optimal pair $((\mun,\rhon,\rgan),\ugan)$, where $\,\ugan\in\uad\,$ and  $\,(\mun,\rhon,\rgan)=
\SALN(\ugan)$. Since $\uad$ is a bounded subset of ${\cal X}$, we may without loss of generality assume
that $\,\ugan\to\ug\,$ weakly-star in ${\cal X}$ for some $\,\ug\in\uad$. 
Then, for the unique solution $(\mu,\rho,\rg, \xi,\xi_\Gamma)$ to \eqref{ss1}--\eqref{ss7} associated
with $\ug$, we conclude from Theorem 3.1 the convergence properties \eqref{conmu}--\eqref{conj2}.
Invoking the optimality of $((\mun,\rhon,\rgan),\ugan)$ for (${\cal P}_{\alpha_n}$),
we then find, for every $v_\Gamma\in\uad$, that 
\begin{align}
\label{tr3.3}
&{\cal J}((\mu,\rho,\rg),\ug)\,=\,{\cal J}({\cal S}_{0}(\ug),\ug)\,\le\,
\liminf_{n\to\infty}\,{\cal J}({\mathcal S}_{\alpha_n}(u^{\alpha_n}_\Gamma),u^{\alpha_n}_\Gamma)
 \nonumber\\[1mm]
&\leq\,\liminf_{n\to\infty}\,{\cal J}({\cal S}_{\alpha_n}(v_{\Gamma}),v_{\Gamma})\, =\,\lim_{n\to\infty} {\cal J}({\cal S}_{\alpha_n}(v_{\Gamma}),v_{\Gamma})\,=\,
{\cal J}({\mathcal S}_{0}(v_{\Gamma}),v_{\Gamma}),
\end{align}  
which yields that $\ug$ is an optimal control for \pier{(${\cal P}_0$)} with the associate state
$(\mu,\rho,\rg,\xi,\xi_\Gamma)$. The assertion is thus proved.
\qed

\vspace{7mm}\quad
Corollary 3.2 does not yield any information on whether every solution to the optimal control problem $({\mathcal{P}}_{0})$ can be approximated by a sequence of solutions to the problems $({\mathcal{P}}_{\alpha})$. 
As already announced in the Introduction, we are not able to prove such a 
general `global' result. Instead, we 
can only give a `local' answer for every individual optimizer of $({\mathcal{P}}_{0})$. For this purpose,
we employ a trick due to Barbu~\cite{Barbu}. To this end, let $\bar u_\Gamma\in\uad$
be an arbitrary optimal control for $({\mathcal{P}}_{0})$, and let $(\bm,\br,\brg,\bar\xi,\bar\xiga)$
be the associated solution quintuple to the state system (\ref{ss1})--(\ref{ss7}) in the sense of 
Theorem 2.1. In particular, $\,(\bm,\br,\brg)={\cal S}_0 (\bar u_\Gamma)$. We associate with this 
optimal control the {\em adapted cost functional}
\begin{equation}
\label{cost2}
\widetilde{\cal J}((\mu,\rho,\rg),\ug):={\cal J}((\mu,\rho,\rg),\ug)\,+\,\frac{1}{2}\,\|u_{\Gamma}-\bar{u}_\Gamma\|^2_{L^2(\Sigma)}
\end{equation}
and a corresponding {\em adapted optimal control problem},

\vspace{4mm}
($\widetilde{\mathcal{P}}_{\alpha}$)\quad Minimize $\,\, \widetilde {\cal J}((\mu,\rho,\rg),\ug)\,\,$
for $\,\ug\in\uad$, subject to the condition that  
\hspace*{12mm}(\ref{als1})--(\ref{als5}) be satisfied.

\vspace{3mm}
With a standard direct argument that needs no repetition here, we can show the following 
result.

\vspace{5mm}
{\sc Lemma 3.3:}\quad\,{\em Suppose that the assumptions} {\sc (A1)}--{\sc (A5)},
{\rm \pier{(\ref{defh})--(\ref{phiat0})} {\em are
satisfied, and let $\alpha\in (0,1]$. Then the optimal control problem} 
$(\widetilde{\cal P}_\alpha)$ {\em admits a solution.}
  
\vspace{5mm}
We are now in the position to give a partial answer to the question raised above. We have the following result.

\vspace{5mm}
{\sc Theorem~3.4:}\,\quad{\em Let the assumptions} {\sc (A1)}--{\sc (A5)}, {\rm \pier{(\ref{defh})--(\ref{phiat0})} {\em be fulfilled, suppose that 
$\bar u_\Gamma\in \uad$ is an arbitrary optimal control of} $({\mathcal{P}}_{0})$ {\em with associated state quintuple $(\bm,\br,\brg,\bar\xi,\bar\xi_\Gamma)$, and let $\,\{\alpha_n\}\subset (0,1]\,$ be
any sequence such that $\,\alpha_n\searrow 0\,$ as $\,n\to\infty$. Then there exist a 
subsequence $\{\alpha_{n_k}\}$ of $\{\alpha_n\}$, and, for every $k\in\nz$, an optimal control
 $\,u_\Gamma^{\alpha_{n_k}}\in \uad\,$ of the adapted problem $(\widetilde{\mathcal{P}}_{\alpha_{n_k}})$
 with associated state triple $(\mu^{\alpha_{n_k}},\rho^{\alpha_{n_k}},
 \rho_\Gamma^{\alpha_{n_k}})$ such that, as $k\to\infty$,}
\begin{align}
\label{tr3.4}
&u_\Gamma^{\alpha_{n_k}}\to \bar u_\Gamma\quad\mbox{strongly in }\,L^2(\Sigma),\\[2mm]
&\mbox{\em the properties {\rm \eqref{conmu}--\eqref{conzig}} are satisfied, where $(\mu,\rho,
\rg,\xi,\xiga)$} \nonumber\\
&\pier{\mbox{\em  is replaced by }\,  (\bar\mu,\br,\brg,\bar\xi,\bar\xi_\Gamma) \,\ \mbox{\em and
the index $\,n\,$ is replaced by $\,n_k$},}\label{pier1}\\[2mm]
\label{tr3.5}
&\widetilde{{\cal J}}((\mu^{\alpha_{n_k}},\rho^{\alpha_{n_k}},\rho_\Gamma^{\alpha_{n_k}}),
u_\Gamma^{\alpha_{n_k}})\to  {\cal J}((\mu,\rho,\rg),\ug)\,.
\end{align}

\vspace{3mm}
{\sc Proof:} \quad\,Let $\alpha_n \searrow 0$ as $n\to\infty$. For any $ n\in\nz$, we pick an optimal control $u_\Gamma^{\alpha_n} \in \uad\,$ for the adapted problem $(\widetilde{\cal P}_{\alpha_n})$ and denote by 
$(\mu^{\alpha_n},\rhon,\rgan)=\SALN(\ugan)$ the associated solution triple of 
problem (\ref{als1})--(\ref{als5}) for $\alpha=\alpha_n$. 
By the boundedness of $\uad$, there is some subsequence $\{\alpha_{n_k}\}$ of $\{\alpha_n\}$ such that
\begin{equation}
\label{ugam}
u_\Gamma^{\alpha_{n_k}}\to \ug\quad\mbox{weakly-star in }\,{\cal X}
\quad\mbox{as }\,k\to\infty,
\end{equation}
with some $\ug\in\uad$. \pier{Thanks to} Theorem 3.1, the convergence properties \eqref{conmu}--\eqref{conzig}
hold true, where $(\mu,\rho,\rg,\xi,\xiga)$ is the unique solution to the state system
\pier{\eqref{ss1}--\eqref{ss7}}. In particular, the pair $(\SO(\ug),\ug)=((\mu,\rho,\rg),\ug)$
is admissible for (${\cal P}_0$). 

\vspace{2mm}\quad
We now aim to prove that $\ug=\bar u_\Gamma$. Once this is shown, then the uniqueness result of 
Theorem 2.1 yields that also $(\mu,\rho,\rg,\xi,\xiga)=(\bm,\br,\brg,\bar\xi,\bar\xi_\Gamma)$, 
which implies that \pier{\eqref{pier1}} holds true. 

\vspace{2mm}\quad Now observe that, owing to the weak sequential lower semicontinuity of 
$\widetilde {\cal J}$, 
and in view of the optimality property
of  $\,((\bm,\br,\brg),\bar u_\Gamma)$ for problem $({\cal P}_0)$,
\begin{align}
\label{tr3.6}
&\liminf_{k\to\infty}\, \widetilde{\cal J}((\mu^{\alpha_{n_k}},\rho^{\alpha_{n_k}},
\rho_\Gamma^{\alpha_{n_k}}),u_\Gamma^{\alpha_{n_k}})
\ge \,{\cal J}((\mu,\rho,\rg),\ug)\,+\,\frac{1}{2}\,
\|u_{\Gamma}-\bar{u}_\Gamma\|^2_{L^2(\Sigma)}\nonumber\\[1mm]
&\geq \, {\cal J}((\bm,\br,\brg),\bar u_\Gamma)\,+\,\frac{1}{2}\,\|u_{\Gamma}-\bar{u}_\Gamma\|^2_{L^2(\Sigma)}\,.
\end{align}
On the other hand, the optimality property of  $\,((\mu^{\alpha_{n_k}},\rho
^{\alpha_{n_k}},\rho_\Gamma^{\alpha_{n_k}}),u^{\alpha_{n_k}}_\Gamma)
\,$ for problem $(\widetilde {\cal P}_{\alpha_{n_k}})$ yields that
for any \pier{$k\in\nz$} we have
\begin{equation}
\label{tr3.7}
\widetilde {\cal J}((\mu^{\alpha_{n_k}},\rho
^{\alpha_{n_k}},\rho_\Gamma^{\alpha_{n_k}}),u^{\alpha_{n_k}}_\Gamma)\, =\,
\widetilde {\cal J}({\cal S}_{\alpha_{n_k}}(u_\Gamma^{\alpha_{n_k}}),
u_\Gamma^{\alpha_{n_k}})\,\le\,\widetilde {\cal J}({\cal S}_{\alpha_{n_k}}
(\bar u_\Gamma),\bar u_\Gamma)\,,
\end{equation}
whence, taking the limit superior as $k\to\infty$ on both sides and invoking (\ref{conj2}) in
Theorem~3.1,
\begin{align}
\label{tr3.8}
&\limsup_{k\to\infty}\,\widetilde {\cal J}((\mu^{\alpha_{n_k}},\rho^{\alpha_{n_k}},
\rho_\Gamma^{\alpha_{n_k}}),u_\Gamma^{\alpha_{n_k}})
\nonumber\\[1mm]
&\le\,\widetilde {\cal J}(\SO(\bar u_\Gamma),\bar u_\Gamma) 
\,=\,\widetilde {\cal J}((\bm,\br,\brg),\bar u_\Gamma)
\,=\,{\cal J}((\bm,\br,\brg),\bar u_\Gamma)\,.
\end{align}
Combining (\ref{tr3.6}) with (\ref{tr3.8}), we have thus shown that 
$\,\frac{1}{2}\,\|u_{\Gamma}-\bar{u}_\Gamma\|^2_{L^2(\Sigma)}=0$\,,
so that $\,\ug=\bar u_\Gamma\,$  and thus also $\,(\mu,\rho,\rg,\xi,\xiga)
=(\bm,\br,\brg,\bar\xi,\bar\xi_\Gamma)$. 
Moreover, (\ref{tr3.6}) and (\ref{tr3.8}) also imply that
\begin{align}
\label{tr3.9}
&{\cal J}((\bm,\br,\brg),\bar u_\Gamma) \, =\,\widetilde{\cal J}((\bm,\br,\brg),\bar u_\Gamma)
\,=\,\liminf_{k\to\infty}\, \widetilde{\cal J}((\mu ^{\alpha_{n_k}},\rho^{\alpha_{n_k}},
\rho_\Gamma^{\alpha_{n_k}}), u_\Gamma^{\alpha_{n_k}})\nonumber\\[1mm]
&\,=\,\limsup_{k\to\infty}\,\widetilde{\cal J}((\mu ^{\alpha_{n_k}},\rho^{\alpha_{n_k}},
\rho_\Gamma^{\alpha_{n_k}}), u_\Gamma^{\alpha_{n_k}}) \,
=\,\lim_{k\to\infty}\, \widetilde{\cal J}((\mu ^{\alpha_{n_k}},\rho^{\alpha_{n_k}},
\rho_\Gamma^{\alpha_{n_k}}), u_\Gamma^{\alpha_{n_k}})\,,
\end{align}                                     
which proves {(\ref{tr3.5})} and, at the same time, also (\ref{tr3.4}). This concludes the proof
of the assertion.\qed

\section{The optimality system}
\setcounter{equation}{0}
In this section, we aim to establish first-order necessary optimality conditions for the optimal control problem $({\mathcal{P}}_{0})$.  This will be achieved by a passage to the limit as $\alpha\searrow 0$ in the first-order necessary optimality conditions for the adapted optimal control problems $(\widetilde{\mathcal{P}}_{\alpha})$ that can by derived as in \cite{CGSneu2} with only
minor and obvious changes. This procedure will yield certain generalized first-order
necessary optimality conditions in the limit. In this entire section, we assume that
$\,h\,$ is given by (\ref{defh}) and that (\ref{phiat0}) and
the general assumptions {\sc (A1)}--{\sc (A5)} are
satisfied. We also assume that a fixed optimal control $\bar u_\Gamma\in \uad$ for 
$({\cal P}_0)$
is given, along
with the corresponding solution quintuple $(\bm,\br,\brg,\bar\xi,\bar\xi_\Gamma)$ of the 
state system
(\ref{ss1})--(\ref{ss7}) established in Theorem 2.1. That is, we have 
$(\bm,\br,\brg)=\SO(\bar u_\Gamma)$ as well as $\bar\xi\in\partial I_{[-1,1]}(\br)$ a.\,e. in $Q$
and $\bar\xi_\Gamma\in \partial I_{[-1,1]}(\brg)$ a.\,e. on $\Sigma$.

\vspace{2mm}\quad
In order to \pier{be} able to take advantage of the analysis performed in \cite[Sect.~4]{CGSneu2},
we impose the following additional compatibility condition:

\vspace{5mm}
(A6) \quad It holds that $\,\,\left(\beta_4(\br(T)-\hat\rho_\oma)\,,\,\beta_5
(\brg(T)-\hat\rho_\Gamma)\right)\in {\cal V}.$

\vspace{5mm}
Obviously, (A6) is fulfilled if $\beta_4=\beta_5$ (especially if $\beta_4=\beta_5=0$) and 
$(\hat\rho_\oma,\hat\rho_\Gamma)\in {\cal V}$. In view of the fact that always 
$(\br(T),\brg(T))\in {\cal V}$, these conditions for the target functions
$\,\hat\rho_\oma\,$ and $\,\hat\rho_\Gamma\,$ seem to be quite reasonable.
 
\vspace{5mm}\quad
 We begin our analysis by formulating the adjoint state system for the adapted
control problem $(\widetilde{\mathcal{P}}_{\alpha})$.
To this end, let us assume that $u_\Gamma^\alpha\in\uad$ is an arbitrary optimal control for 
$(\widetilde{\mathcal{P}}_{\alpha})$ and that $(\mal,\ral,\rgal)$ 
is the solution triple to the associated state system (\ref{als1})--(\ref{als5}). In particular, 
$(\mal,\ral,\rgal)={\cal S}_\alpha(\ugal)$,  the solution has the 
regularity properties (\ref{regmu})--(\ref{regrg}), and it satisfies the global bounds
\eqref{albo1}, \eqref{albo2}, as well as the separation property \eqref{separ}.
Moreover, it follows from \cite[Thm.~4.2]{CGSneu2} that the associated adjoint system
\begin{align}
\label{as1}
&-(1+2g(\ral))\,p_t^\alpha- g'(\ral)\,\rho^\alpha_t\,\pal-\Delta\pal\,=\,g'(\ral)\,\qal
+\beta_1(\mal-\hat\mu_Q)\nonumber \\
&\hskip7cm\mbox{a.\,e. in }\,Q,\\[2mm]
\label{as2}
&\qquad \pn\pal=0\quad\mbox{a.\,e. on }\,\Sigma,\quad \pal(T)=0\quad\mbox{a.\,e. in }\,\oma,\\[2mm]
\label{as3}
&-q^\alpha_t-\Delta\qal+(\vp(\alpha)\,h''(\ral)+\pi'(\ral)-\mal\,g''(\ral))\,\qal
\nonumber\\[1mm]
&\qquad=\,g'(\ral)(\mal\,p_t^\alpha-\mu_t^\alpha\,\pal)+\beta_2(\ral-\hat\rho_Q)
\quad\mbox{a.\,e. in }\,Q,\\[2mm]
\label{as4}
&\pn\qal-\pt\qal-\dega\qgal+(\vp(\alpha)\,h''(\rgal)+\pigs(\rgal))\,\qgal
\,=\,\beta_3(\rgal-\hat\rho_\Sigma),\nonumber\\[1mm]
&\hskip4.5cm\mbox{and }\,\qgal=q^\alpha_{|\Sigma},\quad\mbox{a.\,e. on }\,\Sigma,\\[2mm]
\label{as5}
&\qal(T)=\beta_4(\ral(T)-\hat\rho_\oma)\quad\mbox{a.\,e. in }\,\oma,
\quad\qgal(T)=\beta_5(\rgal(T)-\hat\rho_\Gamma)\nonumber \\
&\hskip7cm\mbox{a.\,e. on }\,\Gamma
\end{align}
has a unique solution $(\pal,\qal,\qgal)$ such that
\begin{align}
\label{eq:4.6}
&\pal\in H^1(0,T;H)\cap C^0([0,T];V)\cap L^2(0,T;W),\\[2mm]
\label{eq:4.7}
&\qal\in H^1(0,T;H)\cap  C^0([0,T];V)\cap L^2(0,T;\hzwei),\\[2mm]
\label{eq:4.8}
&\qgal\in H^1(0,T;\Hg)\cap  C^0([0,T];\Vg)\cap L^2(0,T;\hzweig).
\end{align}

In addition, as in the proof of \cite[Cor.~4.3]{CGSneu2}, it follows the validity of 
the variational inequality 
\begin{equation}
\label{eq:4.9}
\int_0^T\!\!\ginto \bigl(q^\alpha_\Gamma + \beta_6\,u_\Gamma^\alpha
+(u_\Gamma^\alpha-\bar u_\Gamma)\bigr)(v_\Gamma-u_\Gamma^\alpha)\,\dg\,\dt \,\ge\,0 \quad
\forall\,v_\Gamma \in \uad\,.
\end{equation}

\vspace{3mm}\quad
We now prove an a priori estimate that will be fundamental for the derivation of the optimality conditions for $(\mathcal{P}_{0})$. To this end, we introduce some further function spaces. At first, we put
\begin{align}
\label{defY}
&Y:=H^1(0,T; V^*)\cap L^2(0,T;V),\quad Y_\Gamma:=H^1(0,T; V_\Gamma^*)\cap L^2(0,T;\Vg),
\\[2mm]
\label{defW}
&{\mathcal W}\,:= \bigl( H^1(0,T; V^*) \times H^1(0,T; V_\Gamma^*) \bigr)
\cap L^2(0,T;{\cal V}),\\[2mm]
\label{defW0}
&{\mathcal W}_0 \,:\,=\{(\eta,\eta_{\Gamma})\in {\mathcal W} :
(\eta {(0)},\eta_{\Gamma}{(0))}=(0,0)\},
\end{align}
which are Banach spaces when equipped with the natural norm 
of $\,Y\times Y_\Gamma$. Moreover, we have the dense and continuous injections
$\,Y\subset L^2(0,T;V)\subset L^2(Q)\subset L^2(0,T;V^*)\subset Y^*\,$ and 
$\,Y_\Gamma\subset L^2(0,T;\Vg)\subset L^2(\Sigma)\subset L^2(0,T;V^*_\Gamma)
\subset Y^*_\Gamma$, where it is understood that
\begin{align}
\label{dupro1}
\left\langle z,v\right\rangle_Y\,&=\,\int_0^T\!\langle z(t),v(t)\rangle_V\dt \nonumber\\
&\quad\mbox{for all \,$z\in L^2(0,T;V^*)\,$ and \,$v\in  L^2(0,T;V)$,}\\
\label{dupro2}
\left\langle z_\Gamma,v_\Gamma\right\rangle_{Y_\Gamma}
\,&=\,\int_0^T\!\langle z_\Gamma(t),v_\Gamma(t)\rangle_{\Vg}\dt \nonumber\\
&\quad\mbox{for all \,$z_\Gamma\in L^2(0,T;V^*_\Gamma)\,$ and \,$v_\Gamma
\in L^2(0,T;\Vg)$}.
\end{align}
We also note that the embeddings $Y\subset C^0([0,T];H)$ and
$Y_\Gamma\subset C^0([0,T];\Hg)$ are continuous. 
Likewise, we have the dense and continuous embeddings 
$\,
{\cal W}\subset L^2(0,T;{\cal V})$ $\subset L^2(0,T;H\times\Hg)
\subset L^2(0,T;{\cal V}^*)\subset {\cal W}^*$, as well as the
continuous injection
$\,{\cal W}\subset C^0([0,T];H\times\Hg)$, which gives the initial condition encoded in \eqref{defW0}
a proper meaning.
Furthermore, since ${\mathcal W}_0$ is a closed subspace of 
$Y\times Y_\Gamma$, 
we deduce that the elements $F \pier{{}= (z, z_\Gamma)} \in {\mathcal W}_0^* $ are exactly those that are of the form 
\begin{align}
\label{eq:4.15}
\left\langle F,(\eta,\eta_\Gamma)\right\rangle_{{\cal W}_0} \,=\,
\left\langle z, \eta  \right\rangle_Y\,+\,
\left\langle z_\Gamma,\eta_\Gamma \right\rangle_{Y_\Gamma} \quad 
\hbox{ for all $(\eta,\eta_\Gamma)\in {\mathcal W}_0$,}
\end{align}
where $z\in Y^*$ and $z_\Gamma\in Y_\Gamma^*$. 
In particular, for $z\in L^2(0,T; V^*)$ and $z_\Gamma\in L^2(0,T;V^*_\Gamma)$ the 
formulas \eqref{dupro1} and \eqref{dupro2} apply.
Observe that these representation formulas allow us to give a 
proper meaning to statements like 
$$
(z^\alpha , z_\Gamma^\alpha ) \to  (z , z_\Gamma )
\quad \hbox{ weakly in } {\mathcal W}^*_0 .  
$$

\quad In addition to the spaces introduced in \eqref{defY}--\eqref{defW0}, we also define
\begin{equation}
\label{defZ}
{\cal Z}:=\left(L^\infty(\pier{0,T;{}}H)\times L^\infty(\pier{0,T;{}}\Hg)\right)\cap L^2(0,T;{\cal V}),
\end{equation}
which is a Banach spaced when endowed with its natural norm.

\vspace{2mm}\quad
 We have the following
result. 

\vspace{5mm}
{\sc Proposition 4.1:} \quad\,{\em Let the general assumptions} {\sc (A1)}--{\sc (A6)},
 {\rm \pier{(\ref{defh})--(\ref{phiat0})}} {\em be satisfied, and let}
\begin{equation}
\label{eq:4.16}
(\lambda^\alpha,\lambda_\Gamma^\alpha)\,:=\,\left(\varphi(\alpha)\,h''(\ral)
\,q^{\alpha}\,, \, \vp(\alpha)\,h''(\rgal)\,q_\Gamma^{\alpha}\right)
\quad\forall\,\alpha\in (0,1].
\end{equation}
{\em Then there exists
a constant $K_3^*>0$, which depends only on the data of the system and on $R$, such
that for all $\alpha\in (0,1]$ it holds}
\begin{align}
\label{eq:4.17}
&\left\|\pal\right\|_{H^1(0,T;H)\cap C^0([0,T];V)\cap L^2(0,T;W)}\,+\,\max_{0\le t\le T}
\left(\|\qal(t)\|_H+\|\qgal(t)\|_{\Hg}\right)\nonumber\\[2mm]
&+\,
\left\|(q^\alpha,q_\Gamma^\alpha)\right\|_{L^2(0,T;{\cal V})}\,+\,
\left\|(\lambda^{\alpha},\lambda^{\alpha}_\Gamma)\right\|_{{\cal W}^*_0}
\,+\,\left\|\left(\partial_t q^{\alpha},
\partial_t q_\Gamma^{\alpha}\right)\right\|_{{\cal W}^*_0} 
\,\le\,K_3^*\,.
\end{align}

\vspace{3mm}
{\sc Proof:} \quad\,In the following, $C>0$ denote positive constants that
may depend on the data of the system but not on $\alpha\in (0,1]$. We make repeated use
of the global estimates \eqref{albo1} and \eqref{albo2} without further reference.

\vspace{2mm}\quad
First, we add
$\,\pal\,$ on both sides of \eqref{as1}, multiply the result by $\,-p_t^\alpha$,
 and integrate over $\oma\times(t,T]$, where $t\in
[0,T)$. Using the fact that $\,\pal(T)=0$, we obtain the inequality
\begin{align}
\label{past1}
\int_t^T\!\!\!\xinto|p_t^\alpha|^2\dx\ds\,+\,\frac 12\,\|\pal(t)\|_V^2\,\le\,I_1+
I_2+I_3,
\end{align} 
where the quantities $I_j$, $1\le j\le 3$, are specified and estimated below. 
At first, Young's inequality yields that
\begin{align}
\label{past2}
I_1:&=\,-\int_t^T\!\!\!\xinto\left(\pal\,+\,\beta_1(\mal-\hat\mu_Q)\right)\,
p_t^\alpha\dx\ds\nonumber\\
&\le\,\frac 15\int_t^T\!\!\!\xinto|p_t^\alpha|^2\dx\ds\,+\,C\,+\,C\int_t^T\!\!\!\xinto
|\pal|^2\dx\ds\,.
\end{align}
Likewise, \pier{we have that}
\begin{align}
\label{past3}
I_2:\,=\,-\int_t^T\!\!\!\xinto g'(\ral)\,\qal\,\rho_t^\alpha\dx\ds\,\le\,
\frac 15\int_t^T\!\!\!\xinto|\rho_t^\alpha|^2\dx\ds\,+\,C\int_t^T\!\!\!\xinto
|\qal|^2\dx\ds\,.
\end{align} 
Moreover, by also invoking H\"older's inequality and the continuity of the embedding
$V\subset L^4(\oma)$, \pier{we deduce that}
\begin{align}
\label{past4}
I_3:&=\,-\int_t^T\!\!\!\xinto g'(\ral)\,\rho_t^\alpha\,\pal\,p_t^\alpha\dx\ds
\,\le\,C\int_t^T\!\|\rho_t^\alpha(s)\|_4\,\|\pal(s)\|_4\,\|p_t^\alpha(s)\|_2\ds
\nonumber\\
&\le \frac 15\int_t^T\!\!\!\xinto|p_t^\alpha|^2\dx\ds\,+\,C\int_t^T\!
\|\rho_t^\alpha(s)\|_V^2\,\|\pal(s)\|_V^2\ds\,,
\end{align}
where the mapping $\,s\mapsto \|\rho_t^\alpha(s)\|_V^2\,\,$ is bounded in
$L^1(0,T)$ uniformly with respect to $\alpha\in (0,1]$. 

\vspace{2mm}\quad
Next, we multiply \pier{\eqref{as3}} by $\,\qal\,$ and integrate over   
$\oma\times (t,T]$, where $t\in [0,T)$. Taking \pier{\eqref{as4}}
into account, we obtain the identity
\begin{align}
\label{past5}
&\frac 12\left(\|\qal(t)\|_H^2+\|\qgal(t)\|^2_{\Hg}\right)
\,+\,\int_t^T\!\!\!\xinto|\nabla q^\alpha|^2\dx\ds \,+\,\int_t^T\!\!\!\ginto
|\nabla_\Gamma\qgal|^2\dg\ds\nonumber\\
&+\,\int_t^T\!\!\!\xinto\vp(\alpha)\,h''(\ral)\,|\qal|^2\dx\ds\,+\,
\int_t^T\!\!\!\ginto\vp(\alpha)\,h''(\rgal)\,|\qgal|^2\dg\ds\nonumber\\
&=\,\frac 12\left(\|\qal(T)\|_H^2+\|\qgal(T)\|_{\Hg}^2\right)
\,
\nonumber\\
&\quad+\,\int_t^T\!\!\!\xinto\left(\mal\,g''(\ral)-\pi'(\ral)\right)|\qal|^2\dx\ds\,
+\,\int_t^T\!\!\!\xinto\beta_2(\ral-\hat\rho_Q)\,\qal\dx\ds\nonumber\\
&\quad  -\int_t^T\!\!\!\ginto \pigs(\rgal)\,|\qgal|^2\dg\ds
\,+\,\int_t^T\!\!\!\ginto\beta_3(\rgal-\hat\rho_\Sigma)\,\pier{\qgal}\dx\ds\nonumber\\
&\quad +\int_t^T\!\!\!\xinto g'(\ral)\,\mal\,p_t^\alpha\,\qal\dx\ds
\,-\,\int_t^T\!\!\!\xinto g'(\ral)\,\mu_t^\alpha\,\pal\,\qal\dx\ds\,.
\end{align} 

Since $\,\vp(\alpha)\,h''\ge 0$, all summands on the left-hand side are nonnegative.
Moreover, invoking \pier{\eqref{as5}} and Young's inequality, it is readily seen that
the first five summands on the right-hand side are bounded by an expression of the form
\begin{align}
\label{past6}
C\,\Big(1\,+\,\int_t^T\!\!\!\xinto|\qal|^2\dx\ds\,+\,
\int_t^T\!\!\!\ginto|\qgal|^2\dg\ds\Big)\,.
\end{align}

\quad It thus remains to estimate \pier{the} last two summands on the right-hand side,
which we denote by $\,J_1\,$ and $\,J_2$, respectively. By virtue of H\"older's
and Young's inequality, we first have that
\begin{align}
\label{past7}
J_1&\le \,C\int_t^T\|\mal(s)\|_\infty\,\|p_t^\alpha(s)\|_2\,\|\qal(s)\|_2\ds
\nonumber\\
&\le\,\frac 15\int_t^T\!\!\!\xinto|p_t^\alpha|^2\dx\ds\,+\,C
\int_t^T\!\!\!\xinto|\qal|^2\dx\ds\,,
\end{align}
while, also using the continuity of the embedding $V\subset L^4(\oma)$,
\begin{align}
\label{past8}
J_2&\le C\int_t^T\!\|\mu_t^\alpha(s)\|_2\,\|\pal(s)\|_4\,\|\qal \pier{(s)}\|_4\ds\nonumber\\
&\le\,\frac 12\int_t^T\!\|\qal(s)\|_V^2\ds\,+\,C\int_t^T\|\mu_t^\alpha(s)\|_H^2
\,\|\pal(s)\|_V^2\ds\,,
\end{align}
where the mapping $\,\,s\mapsto \|\mu_t^\alpha(s)\|_H^2\,\,$ is known to be bounded 
in $L^1(0,T)$, uniformly in $\alpha\in (0,1]$. Therefore, combining the 
estimates \eqref{past1}--\eqref{past8}, we obtain from Gronwall's lemma, taken
backward in time, the estimate
\begin{align}
\label{past9}
\left\|\pal\right\|_{H^1(0,T;H)}\,+\,\max_{0\le t\le T}\left(
\|\pal(t)\|_V\,+\,\|\qal(t)\|_H\,+\,\|\qgal(t)\|_{\Hg}\right)
 \nonumber \\
 +\,\left\|(\qal,\qgal)\right\|_{L^2(0,T;{\cal V})}\,\le\,C.
\end{align}

\vspace{1mm} Now observe that 
\begin{align*}
&\left\|g'(\ral)\,\rho_t^\alpha\,\pal\right\|^2_{L^2(Q)}\,\le\,
C\txinto|\rho_t^\alpha|^2\,|\pal|^2\dx\dt\\
&\le\,C
\int_0^T\!\!\|\rho_t^\alpha(s)\|_4^2\,\|\pal(s)\|_4^2\ds\,\le\,C.
\end{align*}
Thus, by comparison in \pier{\eqref{as1}}, \pier{we find out that} $\,\|\Delta\pal\|_{L^2(Q)}\,\le\,C$,
whence, by virtue of \pier{\eqref{as2} and} standard elliptic estimates,
\begin{equation}
\label{past10}
\left\|\pal\right\|_{L^2(0,T;W)}\,\le\,C.
\end{equation}

\vspace{2mm}\quad
Next, we derive the bound for the time derivatives. To this end, let 
$\,(\eta,\eta_\Gamma)\in {\mathcal W}_0 \,$ be arbitrary. 
Using the continuity of the embeddings $\,Y \subset C^0([0,T];H)$ and $\,Y_\Gamma\subset
C^0([0,T];\Hg)$, and invoking the estimate 
(\ref{past9}), we obtain from integration by parts  
that
\begin{align}
&\left\langle(\partial_t q^\alpha,\pt \qgal),(\eta,\eta_\Gamma)\right\rangle
_{{\mathcal W}}\,=\,\texinto q_t^\alpha\,\eta\dx\dt\,+\,\teginto\pt q_\Gamma^\alpha
\,\eta_\Gamma\dg\dt\nonumber\\[2mm] 
&=\xinto\!\!\qal(T)\,\eta(T)\dx\,+\ginto \!q^\alpha_\Gamma(T)\,\eta_\Gamma(T)\dg
\nonumber\\[2mm] 
&\quad-\int_0^T\!\!\!\left\langle\eta_t(t),\qal(t)\right\rangle_V\dt
\,-\int_0^T\!\!\!\left\langle\pt\eta_\Gamma(t),\qgal(t)\right\rangle_{\Vg}\dt
\nonumber\\[2mm]
&\le\,\|\qal(T)\|_H\,\|\eta(T)\|_H\,+\,\|\qgal(T)\|_{\Hg}\,\|\eta_\Gamma(T)\|
_{\Hg}\nonumber\\
&\quad+\int_0^T\!\|\eta_t(t)\|_{V^*}\,\|\qal(t)\|_V\dt\,+\int_0^T\!\|\pt\eta_\Gamma(t)\|_{V_\Gamma^*}\,\|\qgal(t)\|_V\dt, 
\nonumber
\end{align}
\pier{whence}
\begin{align}
&\pier{\left\langle(\partial_t q^\alpha,\pt \qgal),(\eta,\eta_\Gamma)\right\rangle_{{\mathcal W}}} 
\nonumber\\[2mm]
&\le C\max_{0\le t\le T}\left(\|\eta(t)\|_H+\|\eta_\Gamma(t)\|_{\Hg}\right)
\nonumber\\
&\quad+C\left\|(\qal,\qgal)\right\|_{L^2(0,T;{\cal V})}
\left(\|\eta_t\|_{L^2(0,T;V^*)}+\|\pt\eta_\Gamma\|_{L^2(0,T;\Vg^*)}
\right)
\le\,C\,\|(\eta,\eta_\Gamma)\|_{{\cal W}_0}\,. \nonumber
\end{align}
We thus have shown that
\begin{equation}
\label{past11}
\left\|
(\partial_t q^\alpha ,\partial_t q_\Gamma^{\alpha}
)\right\|_{{\cal W}^*_0}
\,\le \,C.
\end{equation}

\vspace{2mm}\quad Now, let $(\eta,\eta_\Gamma)\in {\mathcal{W}}_0$ be arbitrary. We define the
functions 
\begin{align}
\label{defval}
&v_1^\alpha:=\left(\mal g''(\ral)-\pi'(\ral)\right)\qal\,+\,g'(\ral)\mal
p_t^\alpha,\quad v_2^\alpha:=-g'(\ral)\mu_t^\alpha\pal,\nonumber\\[1mm]
&w^\alpha:=-\pigs(\rgal)\qgal\,.
\end{align}
Multiplying \eqref{as3} by $\eta$, and invoking \eqref{as4}, we then easily 
infer the identity
\begin{align}
\label{past12}
&\left\langle (\lambda^\alpha,\lambda_\Gamma^\alpha),(\eta,\eta_\Gamma)
\right\rangle_{{\cal W}_0}
\,=\,\texinto\lambda^\alpha\,\eta\dx\dt\,+\,\teginto\lambda_\Gamma^\alpha
\,\eta_\Gamma\dg\dt\nonumber\\
&=\texinto\eta\,q_t^\alpha\dx\dt\,+\teginto\eta_\Gamma\,\pt q_\Gamma^\alpha
\nonumber\\
&\quad-\texinto\nabla q^\alpha\cdot\nabla\eta\dx\dt\,
-\teginto\nabla_\Gamma q_\Gamma^\alpha\cdot\nabla_\Gamma\eta_\Gamma\dg\dt\nonumber\\
&\quad+\texinto v_1^\alpha\,\eta\dx\dt\,+\texinto v_2^\alpha\eta\dx\dt
\,+\teginto w^\alpha\,\eta_\Gamma\dg\dt \nonumber\\
&\quad\pier{+\texinto \beta_2 (\rho^\alpha - \hat\rho_Q)\dx\dt
\,+\teginto \beta_3 (\rho^\alpha_\Gamma - \hat\rho_\Sigma)\dg\dt.}
\end{align}

Now observe that $\,v^\alpha_1\,$ and $\,w^\alpha\,$ are known to bounded in
$L^2(Q)$ and in $L^2(\Sigma)$, respectively, uniformly in $\alpha\in (0,1]$. Also,
using the continuity of the embedding $\hzwei\subset\lio$, \pier{we have that}
\begin{align}
&\texinto v_2^\alpha\,\eta\dx\dt\,\le\,C\int_0^T\!\!\|\mu_t^\alpha(t)\|_2\,
\|\eta(t)\|_2\,\|\pal(t)\|_\infty\dt\nonumber\\[2mm]
&\le\,C\,\max_{0\le t\le T}\,\|\eta(t)\|_H\,\|\mu_t^\alpha\|_{L^2(Q)}\,
\|\pal\|_{L^2(0,T;\hzwei)}\,\le\,C\,\|\eta\|_Y\,.
\end{align}

Therefore, taking \eqref{past9} and \eqref{past11} into account, we have
shown that
\begin{equation}
\label{past13}
\|(\lambda^{\alpha},\lambda_\Gamma^{\alpha})\|_{{\cal W}^*_0}\,\le\,C\,.
\end{equation}
This concludes the proof of the assertion.
\qed

\vspace{7mm}\quad After these preliminaries, we are now in a position
to establish first-order necessary optimality conditions for $({\mathcal{P}}_0)$
by performing a limit as $\alpha\searrow 0$ in the approximating problems. To this end, 
recall that a fixed optimal control $\bar u_\Gamma\in \uad$ for $({\cal P}_0)$, along
with a solution quintuple $(\bm,\br,\brg,\bar\xi,\bar\xiga)$ of the associated 
state system (\ref{ss1})--(\ref{ss7}) is given.

\vspace{2mm}\quad
Now, we choose an arbitrary sequence $\{\alpha_n\}$ such that $\alpha_n
\searrow0$ as $n\to\infty$. By virtue of Theorem 3.4, we can find a subsequence,
which is again indexed by $\,n$, such that, for any $n\in\nz$,
we can find an optimal control $u_\Gamma^{\alpha_n}\in\uad$ for
$(\widetilde{\mathcal{P}}_{\alpha_n})$ with associated state triple 
$(\mun,\rhon,\rgan)$ that satisfies the convergence properties 
(\ref{tr3.4})--(\ref{tr3.5}). From \cite[Sect.~8, Cor.~4]{Simon}, 
\pier{without loss of generality we may} assume that 
\begin{align}
\label{past14}
\mun&\to \bm\quad\mbox{strongly in }\,C^0([0,T];L^p(\oma))\quad\mbox{for  }\,1\le p<6,\\[1mm]
\label{past15}
\rhon&\to\br\quad\mbox{strongly in }\,\CQ,\quad\rgan\to\brg\quad
\mbox{strongly in }\,\CS, 
\end{align}
which \pier{entail} that
\begin{align}
\label{past16}
\Psi(\rhon)&\to\Psi(\br)\quad\mbox{strongly in }\,\CQ\quad\mbox{for }\,
\Psi\in\{g,g',g'',\pi,\pi'\}\\[1mm]
\Psi_\Gamma(\rgan)&\to\Psi_\Gamma(\br)\quad\mbox{strongly in }\,\CS
\quad\mbox{for }\,\Psi\in\{\pig,\pigs\}.
\end{align}

\quad
Moreover, thanks to Proposition 4.1 and to \cite[Sect.~8, Cor.~4]{Simon}, we may assume that the associated
adjoint variables $(p^{\alpha_n},q^{\alpha_n},q_\Gamma^{\alpha_n})$ satisfy
\begin{align}
\label{past17}
&p^{\alpha_n}\to p\quad\mbox{weakly\pier{-star} in }\,H^1(0,T;H)\cap \pier{L^\infty(0,T;V)}
\cap L^2(0,T;W)\nonumber\\
&\hspace*{16mm}\mbox{and strongly in }\,C^0([0,T];L^p(\oma))\quad\mbox{for }\,1\le p<6,\\[2mm]
\label{past18}
&(q^{\alpha_n},q_\Gamma^{\alpha_n})\to (q,q_\Gamma)\quad\mbox{weakly-star in }\,
\pier{{\cal Z},}\\[2mm]
\label{past18bis}
&\pier{(\pt q^{\alpha_n},\pt q_\Gamma^{\alpha_n})\to (\pt q,\pt q_\Gamma)\quad\mbox{weakly in }\,
{\cal W}_0^*\, ,}\\[2mm]
\label{past19}
&(\lambda^{\alpha_n},\lambda_\Gamma^{\alpha_n})\to (\lambda,\lambda_\Gamma)
\quad\mbox{weakly in } {\cal W}_0^*\,,
\end{align}
for suitable limits $\,(p,q,q_\Gamma)\,$ and \,$(\lambda,\lambda_\Gamma)$,
where $\lambda\in Y^*$ and $\lambda_\Gamma\in Y_\Gamma^*$, 
as explained \pier{around \eqref{eq:4.15}}. Obviously, \eqref{past17} implies that $\,\pn p=0\,$ almost everywhere on $\Sigma$ and
$p(T)=0$ almost everywhere in $\oma$. 
Therefore, passing to the limit as
$n\to\infty$ in the variational inequality (\ref{eq:4.9}), written for $\alpha_n$, 
$n\in\nz$, we obtain that $(p,q,q_\Gamma)$ satisfies
\begin{equation}
\label{past20}
\teginto\!(q_\Gamma\,+\,\beta_6\,{\bar u}_\Gamma)\,(v_\Gamma-
{\bar u}_\Gamma)\dg\dt\,\ge\,0 \quad
\forall\,v_\Gamma\in\uad.
\end{equation}

\vspace*{2mm}\quad
Next, we aim to show that in the limit as $n\to \infty$ a limiting adjoint 
system for $({\cal P}_0)$
is satisfied. At first, it easily follows from the convergence properties stated above that
\begin{align}
\label{past21}
&g(\rhon)\,p_t^{\alpha_n}\to g(\br)\,p_t, \quad g'(\rhon)\,\rho_t^{\alpha_n}\,p^{\alpha_n}\to g'(\br)\,\br_t\,p,
\quad g'(\rhon)\,q^{\alpha_n} \to g'(\br)\,q,
\end{align}
all weakly in $L^1(Q)$. It thus follows, by taking the limit as $n\to\infty$ in \eqref{as1} and \eqref{as2},
that the limits $p,q$ satisfy 
\begin{align}
\label{limit1}
-\,(1+2g(\br))\,p_t-g'(\br)\,\br_t\,p-\Delta p\,=\,g'(\br)\,q+\beta_1(\bm-\hat\mu_Q)\quad\mbox{a.\,e. in }\,Q,\\[1mm]
\label{limit2}
\pn p=0\quad\mbox{a.\,e. on }\,\Sigma,\quad p(T)=0\quad\mbox{a.\,e. in }\,\oma.\hspace*{20mm}
\end{align}

\vspace{2mm}\quad
The limiting equation corresponding to \eqref{as3}--\eqref{as5} has to be formulated in a weak form.
 To this end, we multiply \eqref{as3}, written for $\alpha_n$, $n\in\nz$, by an
 arbitrary $(\eta,\eta_\Gamma)\in 
{\cal W}_0\,$   and 
integrate the resulting equation over $Q$. Integrating by parts with respect to 
time and space, and invoking the endpoint conditions for $q$ and $q_\Gamma$, as
well as the zero initial conditions for $(\eta,\eta_\Gamma)$,  we arrive at the identity
\begin{align}
\label{past22}
&\texinto\lambda^{\alpha_n}\,\eta\dx\dt\,+\!\teginto \lambda_\Gamma^{\alpha_n}\,\eta_\Gamma\dg\dt
\,+\!\int_0^T\!\langle {\partial_t} \eta (t),q^{\alpha_n}(t)\rangle_V \dt\nonumber\\[2mm]
&+\int_0^T\!\langle
\partial_t\eta_\Gamma(t),q_\Gamma^{\alpha_n}(t)\rangle_{V_\Gamma}\dt
 +\texinto \nabla q^{\alpha_n}\cdot \nabla\eta\dx\dt
\,+\,\teginto \nabla_\Gamma q_\Gamma^{\alpha_n}\cdot\nabla_\Gamma \eta_\Gamma\dg\dt\nonumber\\[2mm]
&-\texinto v_1^{\alpha_n}\,\eta\dx\dt \,-\texinto v_2^{\alpha_n}\eta\dx\dt\,-\teginto w^{\alpha_n}\dg\dt
\nonumber\\[2mm]
&= \,\beta_2\texinto ({\rho^{\alpha_n}}-\hat\rho_Q)\,\eta\dx\dt\,+\,\beta_3\teginto
(\rho^{\alpha_n}_\Gamma-\hat\rho_\Sigma)\,\eta_\Gamma
\dg\dt\nonumber\\[2mm]
&\quad+ \,\beta_4\xinto(\rhon(T)-\hat\rho_\oma)\,\eta(T)\dx \,+\,\beta_5\ginto(\rgan(T)-\hat\rho_\Gamma)
\,\eta_\Gamma(T)\dg\,.
\end{align}

\vspace{2mm}
Now, owing to \eqref{dupro1}--\eqref{eq:4.15}, the sum of the first two integrals on the left-hand side
of  \pier{\eqref{past22} is equal to} $\,\langle (\lambda^{\alpha_n},\lambda_\Gamma^{\alpha_n}),(\eta,\eta_\Gamma)\rangle
_{{\cal W}_0}$, which, by \eqref{past19}, converges to $\,\langle (\lambda,\lambda_\Gamma),(\eta,\eta_\Gamma)
\rangle_{{\cal W}_0}$. Moreover, it is straightforward to verify \pier{(and this may be left to the reader)} that
also the remaining integrals in \eqref{past22} converge. We therefore obtain, for every 
$ \,(\eta,\eta_\Gamma)\in {\cal W}_0$,
\begin{align}
\label{past23}
&\langle (\lambda,\lambda_\Gamma)(\eta,\eta_\Gamma)\rangle_{{\cal W}_0}
\,+\!\int_0^T\!\langle {\partial_t} \eta (t),q(t)\rangle_V \dt \,+\int_0^T\!\langle
\partial_t\eta_\Gamma(t),\pier{q_\Gamma}(t)\rangle_{V_\Gamma}\dt\nonumber\\[2mm]
& +\texinto \nabla q\cdot \nabla\eta\dx\dt
\,+\,\teginto \nabla_\Gamma q_\Gamma\cdot\nabla_\Gamma \eta_\Gamma\dg\dt
\,+\teginto\pigs(\brg)\,q_\Gamma\,\eta_\Gamma\dg\dt
\nonumber\\[2mm]
&+\texinto\left[(\pi'(\br)-\bm\,g''(\br))\,q\,+\,g'(\br)\,(\bm_t\,p-\bm\,p_t)\right]\,\eta\dx\dt
\nonumber\\[2mm]
&= \,\beta_2\texinto (\br-\hat\rho_Q)\,\eta\dx\dt\,+\,\beta_3\teginto
(\brg-\hat\rho_\Sigma)\,\eta_\Gamma
\dg\dt\nonumber\\[2mm]
&\quad+ \,\beta_4\xinto(\br(T)-\hat\rho_\oma)\,\eta(T)\dx \,+\,\beta_5\ginto(\brg(T)-\hat\rho_\Gamma)
\,\eta_\Gamma(T)\dg\,.
\end{align}

\vspace*{3mm}
Next, we show that the limit pair $\,((\lambda,\lambda_\Gamma),(q,q_\Gamma))\,$ satisfies some sort
of a complementarity slackness condition. To this end, observe that \pier{(cf.~\eqref{eq:4.16})}
for all $n\in\nz$ we obviously have
$$
\texinto\lambda^{\alpha_n}\,q^{\alpha_n}\dx\dt\,=\,\texinto\varphi(\alpha_n)\,h''
(\rhon)
\,|q^{\alpha_n}|^2\dx\dt\,\ge\,0\,.$$
An analogous inequality holds for the corresponding boundary terms. \pier{Hence, it is found that}
\begin{equation}
\label{slack}
{\liminf_{n\to\infty} \texinto\lambda^{\alpha_n}\,q^{\alpha_n}\dx\dt\, \ge 0,\quad 
\liminf_{n\to\infty} \teginto\lambda_\Gamma^{\alpha_n}\,q_\Gamma^{\alpha_n}\dg\dt\,\ge\,0\,.}
\end{equation}
   
\vspace*{3mm}
Finally, we derive a relation which gives some indication that the limit $(\lambda,\lambda_\Gamma)$ 
should somehow be  concentrated
on the set where $\,|\br|=1\,$ and $\,|\brg|=1$ (which, however, we cannot prove rigorously). 
To this end, we test the pair $\,(\lambda^{\alpha_n},\lambda_\Gamma^{\alpha_n})\,$  by the function
$$\,\left((1-(\rhon)^2)\,\phi, (1-(\rgan)^2)\,\phi_\Gamma\right)\,$$ that belongs to
${\cal V}$, since} $\,(\phi,\phi_\Gamma)\,$ is any smooth test function satisfying 
\begin{equation}
\label{eq:4.43}
(\phi(0),\phi_\Gamma(0))
=(0,0), \quad \int_\oma (1-(\rhon)^2)\,\phi(t)\,\dx=0 \quad\forall\,t\in [0,T].
\end{equation}
 As $\,h''(r)= 2 / \left( 1-r^2 \right) $ for every $r\in (-1,1)$, we obtain that
\begin{align}
\label{eq:4.44}
&\lim_{n\to\infty}\left(\texinto\lambda^{\alpha_n}\,(1-(\rhon)^2)\,\phi\dx\dt
\,,\teginto\lambda_\Gamma^{\alpha_n}\,(1-(\rgan)^2)\,\phi_\Gamma
\dg\dt\right)\nonumber\\[2mm]
&=\,\lim_{n\to\infty}\left(2\texinto \varphi(\alpha_n)\,q^{\alpha_n}\,\phi\dx\dt\,,\,
2\teginto\pier{\varphi(\alpha_n)}\, q_\Gamma^{\alpha_n}\,\phi_\Gamma\dg\dt\right)\,=\,(0,{0})\,.\quad
\end{align}

\vspace{5mm}
We now collect the results established above. We have the following statement.

\vspace{5mm}
{\sc Theorem~4.2:}\,\quad{\em Let the assumptions {\sc (A1)}--{\sc (A6)} \pier{and}}
\pier{\eqref{defh}--\eqref{phiat0}} {\em be satisfied.
Moreover, let $\,\bar u_\Gamma\in\uad$ be an optimal control for $({\cal P}_0)$ with \pier{the
associated quintuple $(\bm,\br,\brg,\bar\xi,\bar\xiga)$ solving} the corresponding state system}
 (\ref{ss1})--(\ref{ss7}) {\em in the sense of Theorem 2.1. Moreover, let $\{\alpha_n\}\subset (0,1]$ be
 a sequence with $\,\alpha_n\searrow0\,$ as $\,n\to\infty\,$ such that there are optimal pairs
$\,((\mun,\rhon,\rgan),\ugan)\,$ for the adapted control problem ($\widetilde{\cal P}_{\alpha_n}$) satisfying}
\eqref{tr3.4}--\eqref{tr3.5} {\em (such
sequences  exist by Theorem 3.4) and having the associated adjoint variables $(p^{\alpha_n},
q^{\alpha_n},q_\Gamma^{\alpha_n})$. Then, for any subsequence $\{n_k\}_{k\in\nz}$ of $\nz$, there are a
subsequence $\{n_{k_\ell}\}_{\ell\in\nz}$ and some quintuple $(p,q,q_\Gamma,\lambda,\lambda_\Gamma)$ such that}
\begin{align}
&p\in H^1(0,T;H)\cap C^0([0,T];V)\cap L^2(0,T;\hzwei),\nonumber\\[1mm]
&\pier{(q,q_\Gamma)\in {\cal Z},\quad (\pt q, \pt q_\Gamma) \in {\cal W}_0^*,} \quad (\lambda,\lambda_\Gamma)\in {\cal W}_0^*,
\end{align}
{\em and such that the relations {\rm \eqref{past17}--\eqref{past19}} are valid (where the sequences are indexed
by $n_{k_\ell}$ and the limits are taken as $\ell\to\infty$). Moreover, the variational inequality}
\eqref{past20} {\em and the adjoint state equations} \eqref{limit1}, \eqref{limit2}, {\em and} \eqref{past23} {\em are
satisfied.}

\vspace{5mm}
{\sc Remark~4.3:}\quad\,Unfortunately, we cannot show that the limit quintuple
$$(p,q,q_\Gamma,\lambda,\lambda_\Gamma)$$ 
solving the adjoint problem associated with the optimal pair 
$$((\bm,\br,\brg,\bar\xi,\bar\xiga),\bar u_\Gamma)$$ is unique. 
Therefore, it may well happen that the limits differ for different
subsequences. However, it turns out  that for any such limit 
$(p,q,q_\Gamma,\lambda,\lambda_\Gamma)$
the component $q_\Gamma$ should satisfy the variational inequality (\ref{past20}). 

\section*{Acknowledgments}
\pier{PC gratefully acknowledges some financial support from the MIUR-PRIN Grant 2015PA5MP7 ``Calculus of Variations'', the GNAMPA (Gruppo
Nazionale per l'Analisi Matematica, la Probabilit\`{a} e loro Applicazioni)
of INDAM (Istituto Nazio\-nale di Alta Matematica) and the IMATI -- C.N.R. 
Pavia.}

 
\end{document}